\documentstyle[12pt,amscd,amssymb,epsf]{amsart}
\begin{document}
\title{The SU(3) Casson Invariant for Integral Homology 3-Spheres}
\author{Hans U. Boden and Christopher M. Herald}
\address{Department of Mathematics, Ohio State University, Mansfield, OH, 44906}
\email{boden@@math.ohio-state.edu}
\address{Department of Mathematics, Swarthmore College, Swarthmore, PA, 19081}
\email{cherald1@@swarthmore.edu}
\date{July 14, 1997}
\begin{abstract}
{We derive a gauge theoretic invariant of integral homology 3-spheres which
counts gauge orbits of irreducible, perturbed flat SU(3) connections with sign
given by spectral flow.
To compensate for the dependence of this sum on perturbations, the invariant
includes contributions from the reducible, perturbed flat orbits.
Our formula for the correction term generalizes that given by Walker 
in his extension of Casson's SU(2) invariant to rational homology 3-spheres.}
\end{abstract}
\maketitle

\newcommand{\const}{\mbox{const}}
\newcommand{\lto}{\longrightarrow}
\newcommand{\al}{\alpha}
\newcommand{\be}{\beta}
\newcommand{\ga}{\gamma}
\newcommand{\de}{\delta}
\newcommand{\ep}{\epsilon}
\newcommand{\th}{\theta}
\newcommand{\la}{\lambda}
\newcommand{\om}{\omega}
\newcommand{\Ga}{\Gamma}
\newcommand{\Om}{\Omega}
\newcommand{\ZZ}{{\Bbb Z}}
\newcommand{\RR}{{\Bbb R}}
\newcommand{\QQ}{{\Bbb Q}}
\newcommand{\CC}{{\Bbb C}}
\newcommand{\BH}{{\Bbb H}}
\renewcommand{\AA}{{\cal A}}
\newcommand{\BB}{{\cal B}}
\newcommand{\GG}{{\cal G}}
\newcommand{\HH}{{\cal H}}
\newcommand{\KK}{{\cal K}}
\newcommand{\MM}{{\cal M}}
\newcommand{\OO}{{\cal O}}
\newcommand{\UU}{{\cal U}}
\newcommand{\rr}{{\frak h}}
\newcommand{\rrp}{{{\frak h}^\perp}}
\newcommand{\ad}{\operatorname{ad}}
\newcommand{\SF}{\operatorname{{\it Sf}}}
\newcommand{\Fred}{\operatorname{Fred}}
\newcommand{\SAFred}{\operatorname{SAFred}}
\newcommand{\id}{\operatorname{id}}
\newcommand{\im}{\operatorname{im}}
\newcommand{\hol}{\operatorname{{\it hol}}}
\newcommand{\tr}{\operatorname{\it tr}}
\newcommand{\CS}{\operatorname{{\it CS}}}
\newcommand{\stab}{\operatorname{Stab}}
\newcommand{\hess}{\operatorname{Hess}}
\renewcommand{\hom}{\operatorname{Hom}}
\newcommand{\codim}{\operatorname{codim}}
\newcommand{\sym}{\operatorname{Herm}} 
\newcommand{\Id}{\operatorname{Id}}
\newcommand{\Span}{\operatorname{span}}
\newcommand{\Spec}{\operatorname{Spec}}
\newcommand{\ind}{\operatorname{ind}}
\newcommand{\coker}{\operatorname{coker}}
\newcommand{\ps}{{\cal F}}

\newtheorem{defn}{Definition}[section]
\newtheorem{lem}[defn]{Lemma}
\newtheorem{thm}[defn]{Theorem}
\newtheorem{prop}[defn]{Proposition}
\newtheorem{cor}[defn]{Corollary}
\newtheorem{theorem}{Theorem}

\section{Introduction}
Since its introduction in 1985,  Casson's invariant  \cite{casson, am}
has been the focus of intense study. For example,
it has been shown  that it extends
as a $\QQ$-valued invariant of oriented 3-manifolds which
retains most of the important properties of the original invariant
(for details, see \cite{walker, lescop} and the references contained therein).
Its relevance to gauge theory was recognized by C.~Taubes, who related 
it to
the Euler characteristic for the instanton homology groups 
defined by A.~Floer \cite{taubes, floer}.
Because Casson's invariant is essentially defined as an 
algebraic count of the number of conjugacy classes of irreducible representations
$\varrho:\pi_1 X \lto SU(2),$ 
it is widely believed that there exists a sequence of related
invariants  $\la_{SU(n)}(X)$ which ``count" the number of conjugacy
classes of
irreducible representations $\varrho:\pi_1 X \lto SU(n)$.
One program for realizing these invariants was proposed
by S.~Cappell, R.~Lee, and E.~Miller
in the research announcement \cite{clm}.

The present article establishes the existence of such an 
invariant for the group $SU(3)$ in case $X$ is an integral 
homology 3-sphere. The main difficulty in defining 
$\la_{SU(n)}(X)$ is that 
one must first perturb so that the space of irreducible representations
is cut out transversely, but the resulting (signed) count
will depend on the perturbation used. To obtain a well-defined invariant, 
one must 
devise a correction term involving only the reducible representations
which compensates for this dependence.  

In extending 
Casson's $SU(2)$ invariant to 
rational homology 3-spheres,
K. Walker gave a formula for the correction term using
the symplectic geometry and stratified structure 
of representation varieties associated
to a Heegaard splitting of the 3-manifold \cite{walker}.
Although the situation of $SU(3)$ representations of integral homology
3-spheres is similar to that of $SU(2)$ representations
of rational homology 3-spheres
(because in both cases there is only one stratum of reducibles to worry 
about),
we adopt a different approach and use instead gauge theory.
This means that
we view conjugacy classes of
representations as gauge orbits of
flat connections  via holonomy 
and study
the moduli space
of solutions to the (perturbed) flatness equation
as the critical set of the (perturbed) Chern-Simons functional.
The appropriate interpretation of our arguments
in the $SU(2)$ case would lead to a
gauge-theoretic formula for Walker's invariant
(cf. \cite{mrowka-walker, lee-li}).

We now give a brief outline of the contents of this paper.
The rest of this section 
presents the fundamental
notions of 3-manifold $SU(3)$ gauge theory and describes 
our main result.
Section 2 introduces the perturbations and 
the perturbed flatness equation.
Section 3 is devoted to establishing structure theorems
for the moduli space of perturbed flat connections and for the parameterized
moduli space.
It is important to notice that regularity for the parameterized moduli space
does not imply that it is smooth; it typically has non-manifold points
which we call bifurcation points. These singularities look locally like 
`T' intersections.

Section 4 introduces the spectral flow orientation on 
the moduli spaces. Subsection 4.4
deserves special mention because it contains a comparison
of the orientations on different strata 
of the parameterized moduli space near a bifurcation point.
This is a key ingredient in 
our main result, which is a formula 
for the $SU(3)$ Casson invariant and the statement that it
defines an invariant of integral homology 3-spheres.
All of this is explained in section 5 (cf. Theorem \ref{mainer}).
The final section contains technical results 
concerning the existence of perturbations 
for $SU(3)$ gauge theory.


Both authors
would like to acknowledge generous postdoctoral support
from McMaster University and 
the Max Planck Institute.
C.H.~ is also grateful to Swarthmore College for a research
grant.
Many thanks to Tomasz Mrowka for suggesting this problem
and for kindly sharing his insight on the subject, and 
also to Andrew Nicas, Brian Hall and Thomas Hunter for 
numerous illuminating conversations.

\subsection{SU(3) gauge theory}

Suppose $X$ is a closed, oriented 3-manifold and $P$ is a principal $SU(3)$
bundle over $X.$ For topological reasons, $P$ is trivial. 
Pick a 
trivialization $P \cong X \times SU(3)$ and denote by
$\Om^p(X;su(3))$ the space of smooth
$p$-forms with values in the adjoint bundle $\ad P \cong X \times su(3).$
Let $\AA$ be the space
of smooth connections in $P;$  
$\AA$ is an affine space modeled on $\Om^1(X;su(3)).$
A gauge transformation is a bundle automorphism $g : P \lto P,$
and the group of smooth gauge transformations
$\GG$ can be identified with $C^\infty(X,SU(3)).$ This group
acts on $\AA$
by $g \cdot A = g  A g^{-1} + g d g^{-1}$
with quotient $$\BB= \AA / \GG.$$ 

As usual, the gauge group action is not free.  
Let $\AA^{*}$ denote the subset of irreducible connections, i.e., 
those with stabilizer $Z(SU(3)) \cong \ZZ_3$, and 
set $\BB^{*}=\AA^{*}/\GG$.  
While $\BB$ 
is singular at gauge orbits  with stabilizer different from $\ZZ_3,$
if $\AA$ and $\GG$ are given the $L^{2}_{1}$ and $L^2_2$
topologies, respectively,
then $\BB^{*}$ inherits the structure of a pre-Banach manifold.  
For the most part, we will omit the 
references to the Sobolev completions in this paper because a detailed 
account of the analysis can be found in \cite{taubes}.

Assume from now on that
$X$ is an integral homology 3-sphere unless otherwise specified.  
Then the stabilizer of any flat connection is isomorphic to $SU(3)$, 
$U(1)$, or $\ZZ_3$ (among nonflat connections, there are two
other possibilities, $U(1)\times  U(1)$ and $S(U(2)\times 
U(1))$).  Let $\AA^{r}$ denote the space of all 
connections 
with stabilizer isomorphic to $U(1)$;  these are the nonabelian
connections which reduce to 
$S(U(2)\times U(1))$ connections.  We adopt the convenient, if not  standard, 
terminology 
whereby $A$ {\em reducible} means $A \in \AA^r.$ 

The quotient $\BB^{r}= \AA^{r}/\GG$, while a singular stratum 
of $\BB$, is itself a smooth manifold.  This may be seen by 
noticing that 
$A\in \AA^{r}$ if and only if it is gauge 
equivalent to a connection whose 1-form 
takes values in $s(u(2)\times u(1))$, and that this 1-form is unique up to gauge 
transformations $g\in C^{\infty} (X,S(U(2)\times U(1)))$.  
Thus $\BB^{r} \cong \AA^*_{S(U(2)\times U(1))} \, / \, 
\GG_{S(U(2)\times U(1))}$.

For $A \in \AA,$ the curvature is the element $F(A) \in \Om^2(X;su(3))$
defined by 
$$F(A) = d A + A \wedge A.$$  
Then $A \in \AA$ is flat in case $F(A) = 0,$ and 
the moduli space of flat connections is
$$\MM = \{ A \in \AA \mid F(A) = 0 \} / \GG \; \subset \, \BB.$$
Set $\MM^*=\MM \cap \BB^{*}$ and $\MM^r=\MM \cap \BB^{r}.$ 
A well known theorem identifies $\MM$ with 
the space of representations
$\varrho:\pi_1 X \lto SU(3)$ modulo conjugation.

The Chern-Simons functional $\CS(A)$ is 
defined by
$$\CS(A) = \frac{1}{8 \pi^2} 
\int_X \tr(A \wedge dA + \tfrac{2}{3} A \wedge A \wedge A).$$
There is an isomorphism
$\pi_0 \GG \cong \ZZ$ given by $g \mapsto \deg g$
(see Proposition \ref{degree}).
If $g \in \GG,$ then $ \CS(g \cdot A) = \deg g + \CS(A),$ thus
$\CS$ descends to a map
$$\CS:\BB \lto \RR / \ZZ = S^1.$$

Choose an orientation and
a Riemannian metric on $X$.  This provides a Hodge star 
operator $* : \Om^p(X;su(3)) \lto \Om^{3-p}(X;su(3))$ and an $L^2$ 
Riemannian metric
on $\AA,$ given by $\langle a, b \rangle_{L^{2}} = -\int_{X}\tr 
(a\wedge *b)$.  
Taking the gradient of $\CS$ with respect to this metric,
one computes that 
$$\nabla \CS(A) = -\frac 1 {4 \pi^2} * F(A),$$
and hence the set of critical points of $\CS,$ modulo $\GG,$
is exactly the moduli space of flat connections $\MM$.

The linearization of the flatness equation $*F(A)=0$ 
is given by the operator $*d_A : \Om^1(X;su(3)) \lto \Om^1(X;su(3)).$
As in \cite{taubes}, we extend this to  the self-adjoint,
elliptic operator
\begin{eqnarray*}
&K_A: \Om^0(X;su(3)) \oplus \Om^1(X;su(3)) \lto \Om^0(X;su(3)) \oplus \Om^1(X;su(3))&\\
&K_A(\xi,a) = (d_A^* a, d_A \xi + * d_A a).&
\end{eqnarray*}
Notice that 
$\ker K_A = \HH^0_A(X;su(3)) \oplus \HH^1_A(X;su(3)),$
the space of $d_A$-harmonic (0+1)-forms.

For $X$ any closed 3-manifold, the moduli space of flat
$SU(3)$ connections $\MM$ is compact and
has expected dimension zero 
since $K_A$ is self-adjoint. 
Achieving transversality requires the
use of perturbations, and we employ 
the same techniques here that were successful in the $SU(2)$ setting
\cite{taubes, herald1,herald2}.

We define a class of admissible perturbation functions in Section 2 by which to vary the 
Chern-Simons functional.  
The construction of an admissible function $h$ 
involves taking a sum of  invariant functions applied to the holonomy 
around a collection of loops (integrated over normal disks of tubular 
neighborhoods of the loops).  
The perturbed Chern-Simons functional is then $\CS_h(A) = \CS(A) + h(A),$
 and a connection is called $h$-perturbed flat if it is a critical point of $\CS_h.$
We show in Section 3 that it is possible to choose an admissible 
function $h$ such that $\MM^{*}_{h}$ and $\MM^{r}_{h}$ are compact 0-dimensional 
submanifolds of $\BB^{*}$  and $\BB^{r}$ consisting of orbits 
that meet a cohomological regularity condition.  

\subsection{Main result}
We begin by recalling from \cite{taubes}
the gauge-theoretic definition of Casson's invariant
$\la(X)$
in case $X$ is an integral homology 3-sphere.
First, choose a small perturbation $h$ 
so that the perturbed flat $SU(2)$ moduli space is a compact, smooth,
oriented 0-manifold. Then the number
of irreducible, perturbed flat connections counted with sign
is seen to be independent of the choice of perturbation $h.$
This follows from the classification of 1-manifolds  
once it is verified that
for generic, one-parameter families 
of perturbations, the irreducible part
of the parameterized $SU(2)$ moduli space 
is a smooth cobordism between the two moduli spaces at either end.
Taubes identified the resulting invariant as $-2$ times Casson's invariant,
normalized as in \cite{am} (see \cite{kkr} for an explanation of the minus sign).

In the $SU(3)$ case, for generic one-parameter families $\rho(t)=h_t$
of perturbations, the irreducible part of the parameterized moduli
space $W^*_\rho$ is an oriented 1-manifold, but it is not generally 
compact.  The reducible part, $W^r_\rho$, is a compact 1-manifold, 
and the union $W^{*}_{\rho} \cup W^{r}_{\rho}$ is compact 
but not smooth.
The problem is illustrated in Figure 1,
where $\rho(t)$ is defined for $t \in [-1,1]$.  
The solid curves depict 
$W_{\rho}^{*}$ and  the dotted curves $W^r_{\rho}$.  
Because of the noncompact ends of $W_{\rho}^{*}$, the parameterized 
moduli space subfails to give a smooth cobordism between $\MM^{*}_{\rho(-1)}$ 
and $\MM^{*}_{\rho(1)}$.  Thus the algebraic sum of perturbed flat 
irreducible orbits is seen to depend on the perturbation in this case.

\begin{figure}[hbt]
\begin{center}
\leavevmode\hbox{%
\epsfxsize=4in
\epsffile{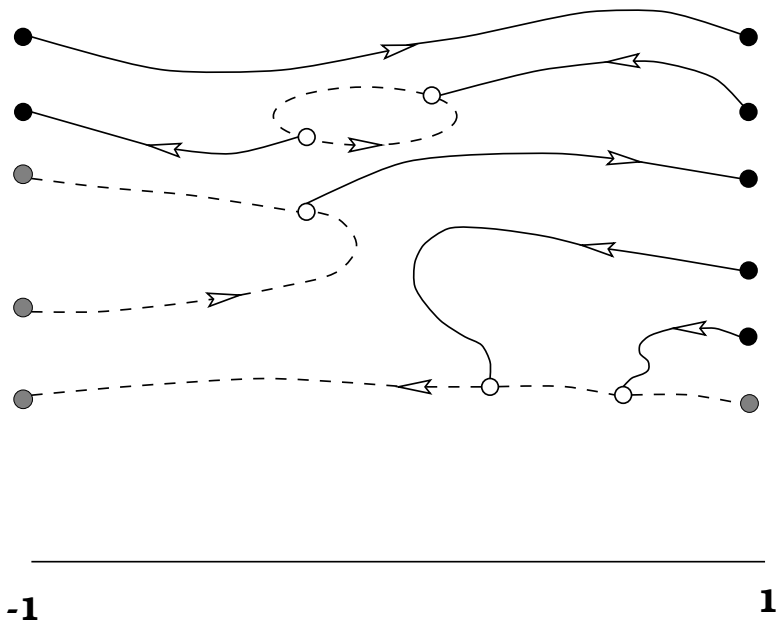}}
\end{center}
\caption{The parameterized moduli space $W^*_\rho \cup W^r_\rho$
projecting vertically to $[-1,1].$ }
\label{figure1}
\end{figure}

The  compactification
$\overline{W^*_\rho}$ 
is obtained by adding certain reducible orbits, called bifurcation points, 
 to the non-compact ends of $W^*_\rho$. 
In Figure \ref{figure1}, the bifurcation points are where the dotted and solid curves meet.
To make the invariant independent of $h,$
one needs a correction term which changes, when the perturbation is 
varied, by the number of  bifurcation points 
on $W^r_\rho$, counted  with sign given by their orientation as 
boundary points of $\overline{W^*_\rho}$.

The oriented spectral flow along $W^r_\rho$ provides
a means to calculate this number, as we now explain.
Let $\rr = s(u(2)\times u(1))$ 
be the 
Lie subalgebra of $su(3)$ and 
$\rrp$ its orthogonal complement, which can be identified 
with $\CC^{2}$.   
For any  reducible connection $A$, the connection 1-form
can be gauge transformed to take values in $\rr$.  
If $A$ is $h$-perturbed flat,
then $\Om^{1}(X;su(3)) = \Om^{1}(X;\rr) \oplus 
\Om^{1}(X;\rrp)$ is the splitting of $T_{A}\AA$ into tangent 
vectors tangent to and normal to 
the reducible stratum. 
For generic paths $\rho,$
the bifurcation points are characterized geometrically
as those reducible orbits in $W_\rho$ where the kernel
of the restriction
of  $K(A,h)$  to the $\rrp$-valued forms jumps up in dimension. 
Such a jump
occurs each time the deformation complex detects a 
tangent vector normal to the reducible stratum.
Hence, in a neighborhood of the bifurcation point in
$W^r_\rho,$ there is a path of eigenvalues of $K(A,h)$
(on $\rrp$-valued forms) crossing zero transversally,
and the sign of its first derivative 
(relative to the orientation on $W_{\rho}^{r}$) coincides with the
boundary orientation of the bifurcation point. 
Note that $\stab A \cong U(1)$ equivariance of $K(A,h)$  forces
the eigenvalue to have multiplicity two.

Choosing the product connection $\th$ as
a reference point for computing all spectral flows, we obtain:
\begin{theorem} \label{mainer}
Suppose $X$ is an integral homology 3-sphere.
For generic small perturbations $h,$
$\MM^*_h$ and $\MM^r_h$ are smooth, compact 
0-manifolds.
Choose representatives $A$ for each orbit $[A] \in \MM_h,$ and
in case $[A] \in \MM^r_h,$ choose also
a flat connection $\widehat{A}$ close to $A$. 
Define $\la_{SU(3)}(X,h)$ to be equal to
$$\sum_{[A] \in \MM^*_h} (-1)^{\SF(\th, A)}
- \frac{1}{2}\sum_{[A] \in \MM^{r}_h}(-1)^{\SF(\th, A)}
(\SF_{\rrp}(\th, A) - 4 \CS(\widehat{A}) + 2),$$
where $\SF$ and $\SF_{\rrp}$ refer to the spectral flow 
of the operator $K(A,h)$ on $su(3)$ and $\rrp$ bundle-valued forms,
respectively. Then for $h$ sufficiently small,
this quantity is independent of $h$ and the Riemannian metric on $X$, 
and gives a well-defined invariant of
integral homology 3-spheres.
\end{theorem}

\noindent
{\it Remark.} 
This theorem will follow from \ref{regularity} and the results in
section 5.

The second sum is our formula for the
correction term. Both
$\SF_{\rrp}(\th, A)$ and $\CS(\widehat{A})$
depend on the choice of representative $A.$ It is only
the difference 
$\SF_{\rrp}(\th, A)  - 4 \CS(\widehat{A})$
which is well-defined on the gauge orbit $[A]$.
The last term in the second sum does not affect the argument
that $\la_{SU(3)}$ is well-defined;
it simply adds 
a certain multiple of the $SU(2)$ Casson invariant to get
a desirable choice of normalization.

As an invariant, $\la_{SU(3)}$ is insensitive to the orientation on $X$.
In general, if $\la_{SU(3)}(X) \neq 0,$ then
$\pi_1 X$ 
admits a non-trivial representation into $SU(2)$ or $SU(3).$
The conjectured rationality of $\CS(\widehat{A})$ would of course
imply that
$\la_{SU(3)}(X) \in \QQ$ as well.

There are many interesting questions raised
by Theorem \ref{mainer}.
The most intriguing is 
what sort of surgery relations (if any) 
does this new invariant satisfy. 
A related question:\footnote{
We are grateful to S. Garoufalidis
for pointing out the connection here.}
is $\la_{SU(3)}$ a finite type 
invariant \cite{oh,gar}? 
By  \cite{murakami},
the Casson-Walker invariant equals 6 times
$\la_1,$ the first Ohtsuki invariant \cite{oh2},
so one is especially interested in any relationship between $\la_{SU(3)}$ and
$\la_2,$ the second Ohtsuki invariant.
Positive results would be interesting for two reasons:
(i) they would render $\la_{SU(3)}$ computable by algebraic means, and (ii)
they would clarify what geometric information the finite type invariants
carry.

There is, of course, still the problem of defining 
the generalized
Casson $SU(n)$ invariants for $n>3.$ A related
problem is 
to extend $\la_{SU(3)}$
to rational homology 3-spheres.
In a different direction, one can attempt to define
$SU(3)$ Floer theory.
We leave these questions to future investigations.
 
\newpage
\section{Perturbations}
In this section, we present the functions that will be used to perturb
the Chern-Simons functional.
After defining the perturbations and
characterizing the perturbed flat connections,
we derive those properties of the first and second
derivative of the perturbation functions which are used later
to prove that the critical set of the perturbed Chern-Simons
functional satisfies certain transversality conditions.  

\subsection{Admissible functions}

This subsection introduces the admissible functions,
which are gauge invariant functions $\AA \lto \RR$
obtained by applying invariant functions $SU(3) \lto \RR$ to 
the holonomy around a collection of loops in $X$. 
We first describe the construction for a single loop.

Each smoothly embedded
based curve $\ell:S^1 \lto X$ defines a holonomy map
$$\hol_{\ell}:\AA \lto SU(3).$$  We can obtain from this a gauge
invariant function $f:\AA \lto \RR$ by composing with an invariant function
$\tau:SU(3)\lto \RR$.  For analytical reasons, it is necessary to 
mollify this function by integrating
against a cut-off function on the 2-disks normal to $\ell$ as follows.

Let $x=(x_1, x_2)$ be coordinates on $D^2$, the unit 2-dimensional disk.
Fix once and for all a radially symmetric 2-form $\eta$  on $D^2$
which vanishes near the boundary and
satisfies $\int_{D^2} \eta =1$. A tubular neighborhood of $\ell$
is an embedded
solid torus $\ga:S^{1} \times D^{2}\lto X.$
For each $x\in D^{2}$, let
$\hol_{\ga} (x,A)$ be the holonomy of $A$ once around
the closed curve $\ga(S^1\times\{x\})$.
For any smooth invariant function $\tau:SU(3) \lto \RR,$ define the
gauge invariant function $p(\ga,\tau):\AA \lto \RR$ by
\begin{equation}
p(\ga,\tau)(A)=\int_{D^2} \tau (\hol_{\ga}(x,A)) \eta(x) dx.
\label{elementary function}
\end{equation}

\begin{defn}\label{admissible functions def}
Fix  $\Ga = \{ \ga_1, \ldots, \ga_n\},$ a set of
embeddings of the solid torus into $X.$
Then an {\bf admissible function relative to}
$\Ga$ is a function $h:\AA \lto \RR$ 
defined by 
$$h(A) = \sum_{i=1}^n p(\ga_i,\tau_i) =
\sum_{i=1}^n \int_{D^2} \tau_i(\hol_{\ga_i}(x,A))\eta(x) dx.$$
where $\tau_i: SU(3) \lto \RR$ is
an invariant 
function of the form $\tau_i=h_i \circ \tr$ for a $C^3$ function
$h_i:\CC\lto \RR$.
Given $\Ga,$ we denote the space of admissible functions by 
$\ps_\Ga$ and note the identification
$\ps_\Ga \cong C^3(\CC,\RR)^{\times n}$ given by
$h \mapsto (h_1,\ldots, h_n).$
For $h \in \ps_\Ga,$ define $\|h\|_{C^3} = \sum_{i=1}^n \|h_i\|_{C^3}.$
\end{defn}
\noindent
 
There is no real loss of generality in considering only
the invariant functions of the type used in the previous definition.
One can see this by the following result, which we have included for
motivation.
\begin{prop}  \label{trace is all}
\smallskip \noindent
\begin{itemize}
\item [(i)]   $\tr : SU(3) \lto \CC$ descends to a one-to-one
map on conjugacy classes.
\item [(ii)]  Any smooth invariant function $\tau:SU(3) \to \RR$
can be written as $\tau = f \circ \tr$ for some
smooth  function  $f : \CC \lto \RR.$
\end{itemize}
\end{prop}  
 
\begin{pf} The characteristic polynomial of $M\in SU(3)$ is given by
$$p_M(\la)= \la^{3} - \tr (M) \la^{2} + \overline{ \tr (M)}
\la -1.$$
Since every matrix in $SU(3)$ is diagonalizable, any two
are conjugate if and only if their eigenvalues coincide,
and (i) follows. 
 
Part (ii) follows from invariant theory.
Consider the case of smooth invariant
functions on $U(3).$ Restricting to a maximal torus $T^3,$
these can be viewed as $S_3$ invariant functions on $T^3,$
where $S_3$ acts 
by permutation of the coordinates.
The inclusion $T^3 \subset \CC^3$ is an equivariant embedding,
and a classical result states that the algebra of invariant
polynomials $P(\CC^n)^{S_n}$ is generated by the elementary,
symmetric functions  $\sigma_1, \ldots, \sigma_n$
(see  Chapter 2A, \cite{weyl}).
This, and Theorem 2 of \cite{schwartz}, proves (ii), since
the $\sigma_i$ are just the coefficients
of the characteristic polynomial, which, for $M \in SU(3),$
are given by $\tr(M)$ and $\overline{\tr(M)}$.
\end{pf} 

\subsection{Perturbed flat connections}
\label{p flat}
In this subsection, we introduce the perturbed flatness equation
and the deformation complex of the perturbed flat moduli space.
Suppose that $\Ga = \{ \ga_{1}, \ldots, \ga_{m} \}$ is
a set of embeddings of the solid torus
into $X.$ All the admissible functions in this section are to
be regarded as admissible relative to $\Ga.$

Pick a Riemannian metric on $X$ and
let $*:\Om^{p}(X;su(3))\to \Om^{3-p}(X;su(3))$
be the Hodge star operator. This defines
an $L^{2}$ inner
product on bundle-valued $p$-forms by
$$\langle \al, \be \rangle_{L^{2}}=-\int_{X}\tr(\al \wedge
*\be )$$
and induces an $L^2$ metric on $\AA.$
For any admissible function
$h:\AA \lto \RR,$ let $\nabla h$ be the gradient of
$h$ with respect to the $L^2$
metric and define
$$\zeta_h: \AA \lto \Om^{1}(X;su(3))$$
by $\zeta_{h}(A) = *F(A) - 4 \pi^2 \nabla h(A)$.
Notice that $\zeta_h(A)$ is just $-4 \pi^2$ times the gradient of the function
from $\AA$ to $\RR$ given by $A \mapsto \CS(A)+h(A).$

\begin{defn}
Suppose $h$ is an admissible function. Then $A \in \AA$ is called
$h$-{\bf perturbed flat} if it satisfies
$$ *F(A) - 4 \pi^2 \nabla h(A) = 0.$$
The {\bf perturbed flat moduli space} is the set
of gauge orbits of perturbed flat connections, i.e.,
$$\MM_h = \zeta_h^{-1}(0)/ \GG.$$
Set
$\MM^{*}_{h}=\MM_{h}\cap \BB^{*}$ and $\MM^{r}_{h}=\MM_{h}\cap \BB^{r}.$
\end{defn}

\begin{defn}\label{param mod}
Suppose $\rho(t), \, -1 \leq t \leq 1,$ 
is a one-parameter family of admissible
functions.
Then the
{\bf parameterized moduli space} is defined as the quotient
$$W_{\rho}= \{(A,t) \in \AA \times [-1,1] \mid \zeta_{\rho(t)}(A)=0 \}/\GG \;
\subset \; \BB \times [-1,1],$$
with slice at $t\in [-1,1]$ given by $\MM_{\rho(t)}\times\{ t\}=
W_{\rho} \cap (\BB \times \{t\}).$
Set $W_{\rho}^{*}=W_{\rho}\cap \left( \BB^{*}\times [-1,1]\right)$ and
$W_{\rho}^{r}=W_{\rho}\cap \left(\BB^{r}\times [-1,1]\right).$
\end{defn}

Since $X$ is an integral homology 3-sphere, 
any reducible {\it flat} connection can be regarded as an irreducible,
flat $SU(2)$ connection. This is no longer true for
perturbed flat reducible connections because they typically
have holonomy in a subgroup conjugate to 
$S(U(2) \times U(1))$
and do not reduce any further.

The linearization  of $\zeta_h$ is given by
$$*d_{A,h} = *d_A - 4 \pi^2 \hess h(A): \Om^1(X;su(3)) \lto \Om^1(X;su(3)).$$
This motivates the final definition of this subsection.

\begin{defn} \label{pert cohom}
Suppose that $h$ is an admissible function
and that $A$ is $h$-perturbed flat.
The {\bf deformation complex} is the 
elliptic Fredholm complex
\begin{eqnarray}\label{complex}
\Om ^{0}(X;su(3)) \stackrel{d_{A}}{\lto}
\Om ^{1}(X;su(3)) \stackrel{*d_{A,h}}{\lto}
\Om ^{1}(X;su(3)) \stackrel{d_{A}^{*}}{\lto}
\Om ^{0}(X;su(3)),
\end{eqnarray}
where $d^*_A$ 
is the $L^2$-adjoint of $d_A.$
The first two  cohomology groups of this complex are
$H^0_A(X;su(3)) = \ker d_A$
 and
$H^{1}_{A,h}(X;su(3)) = \ker *d_{A,h} / \im d_A.$
Notice that this is a self-adjoint complex, and so 
cohomological groups of complementary dimensions are identified.
\end{defn}

Of course, if $h=0$, then (\ref{complex}) is
just the twisted de Rham complex
with the second half rewritten using duality.
We will represent $H^0_A(X;su(3))$ and $H^1_{A,h}(X;su(3))$
by the spaces  
$\HH^0_A(X;su(3))$ and $\HH^1_{A,h}(X;su(3))$
of harmonic forms,
where a 1-form $a$ is {\em harmonic} if
$d_{A} a = 0$ and $*d_{A,h}(a) =0$.
Geometrically, the former cohomology group is the Lie algebra
of $\stab(A),$ while the latter is
the kernel of the linearized
perturbed flatness equation restricted 
to the tangent space to the slice of the gauge group action.

Given a complex line $V \subset \CC^{3}$, we can
decompose $\CC^{3}$ into $V$ and $V^{\perp}.$  This gives an
identification, typically different from the standard one, between
$\CC^{3}$ and $\CC \oplus \CC^2$.
This engenders a corresponding decomposition of 
the Lie algebra as $su(3)= \rr \oplus \rr^\perp$, isomorphic 
(as a vector space) to 
$s(u(2)\times u(1)) \oplus \CC^{2}$. 
For example, for the standard decomposition,
$$
\rr = \left\{ \left(\begin{array}{ccc}  
i(a+b) & c+id & 0\\
                          -c + id & i(a -b)& 0\\
                           0& 0 & -2ia 
\end{array}\right) \right\} 
                           \mbox{ and } \rrp = 
                           \left\{ \left(\begin{array}{ccc}  0 & 0 & 
                           z_{1}\\
                          0 & 0& z_{2}\\
                           -\overline{z}_{1}& -\overline{z}_{2}& 0 
                           \end{array}\right) \right\} 
 .$$
In general, $\rr$ and $\rr^\perp$ are given by conjugating
the above subspaces.

If $A$ is a connection in the bundle $P = X\times SU(3)$  and
$\stab A \cong U(1),$ 
then the action of $\stab A$ on the canonical $\CC^3$ bundle
$E \lto X$
decomposes each fiber of $\ad P$
in a similar manner.  We  shall use 
the notation $\rr$ and $\rrp$  without indicating
the actual dependence of the splitting of $\ad P$ on 
the subgroup $\stab A \subset \GG$;
one can always
gauge transform $A$ into $\AA_{S(U(2)\times U(1))}$
and then $\stab A$ would just give the standard decomposition.

For $A \in \AA^r,$ we decompose 1-forms in a similar manner, and
$\Om^{1}(X;su(3))=\Om^{1}(X;\rr) \oplus \Om^{1}(X;\rrp)$
is a geometric splitting of the 
tangent space
$T_{A}\AA$ into vectors tangent to
the reducible stratum $\AA^{r}$ and vectors normal to that stratum.
If $A$ is $h$-perturbed flat, this leads to
a splitting of the
cohomology groups as
$$\HH^{*}_{A,h}(X;su(3)) = \HH^{*}_{A,h}(X; \rr )
\oplus\HH^{*}_{A,h}(X;\rrp).$$

For convenience, set $\Om^{0+1}(X;su(3)) = \Om ^{0}(X;su(3))\oplus \Om ^{1}(X;su(3)).$ 
We can fold 
the deformation complex (\ref{complex})
up into a single operator
$$K(A,h):\Om ^{0+1}(X;su(3))
\lto \Om ^{0+1}(X;su(3))$$
by setting, for $(\xi, a) \in \Om ^{0}(X;su(3))\oplus \Om ^{1}(X;su(3)),$
$$K(A,h)(\xi, a) = (d_{A}^{*} a, d_{A}\xi + *d_{A,h}(a)).$$
Notice that $K(A,h)$ is a self-adjoint elliptic operator (with 
appropriate Sobolev norms on the domain and range).  
When $A$ is reducible, the operator $K(A,h)$
respects the decomposition of $\Om^{1}(X;su(3))$ described above. 
In particular, in Sections
\ref{orientations and spectral flow} and \ref{the invariant},
we use this to split the spectral flow of $K(A,h)$.

\subsection{The calculus of admissible functions}
In this subsection, we describe the first and second derivatives of
functions $f:\AA \lto \RR$
obtained by composing the holonomy around a loop with an invariant function
$\tau:SU(3) \lto \RR$ as in eqn.~(\ref{elementary function}).

For such functions, these computations can all be performed on
the pullback bundles over $S^1.$
Hence, throughout this section, $\AA$ denotes the
space of connections on the bundle $P= S^1\times  SU(3) .$
Parameterize the circle by
$f:[0,1] \to S^{1}$,  $f(u)=e^{2\pi i u }.$
For $A\in \AA $,
let $\hol (A) \in SU(3)$ be the holonomy once around the circle
in a counterclockwise direction,
based at $1=f(0)$.  

The derivatives
of $\hol (A) $ may be computed
as follows.
For $A\in \AA,$
parallel translation by $A$ defines a trivialization of
the pullback bundle $f^*(\ad P),$ which
identifies tangent vectors in $T_A\AA$  with functions
$a:[0,1]\to su(3)$.

\begin{prop}  \label{derivative of holonomy}
Suppose $A\in \AA $ and $a, b\in T_{A}\AA.$ Then
\begin{itemize}
\item [(i)] $
\left. \frac{d}{dt} \hol(A+ta)\right|_{t=0}=\hol (A) \int_{0}^{1} a(\nu) d\nu
,
$
\item [(ii)]
$\left.\frac{\partial^{2}}{\partial s \partial t}
\hol(A+ta+sb)\right|_{(0,0)}
= \hol(A) \int_{0}^{1}
\int _{0}^{\nu} (a(\nu) b(\mu)+ b(\nu) a(\mu)) d\mu d\nu.  $
\end{itemize}
\end{prop}

\begin{pf}
We prove (ii) and leave (i)  as
an exercise for the reader.

Let $P(s,t;u)\in SU(3)$ denote the parallel translation with respect to
the fixed trivialization from $0$ to
$u$ along the interval by the connection $A + sa + tb$.
Then $P(s,t;u)$ satisfies the differential equation
\begin{equation} \label{curly}
\tfrac{\partial }{\partial u} P(s,t;u) + (sa(u) +
tb(u))P(s,t;u) =0.
\end{equation}
Applying $\frac{\partial^2}{\partial s \partial t}$ to (\ref{curly})
at $(s,t)=(0,0),$ we obtain
$$
\left.{\tfrac{\partial}{\partial u}}\left(
{\tfrac{\partial^2}{\partial s \partial t}}
P(s,t;u)\right|_{(0,0)}\right)
+\left. a(u) {\tfrac{\partial}{\partial t}} P(0,t;u)\right|_{t=0}
+  b(u) \left. {\tfrac{\partial}{\partial s}} P(s,0;u)\right|_{s=0}=0.$$
Integrating with respect to $u$ and commuting mixed partials gives
$$\left.{\tfrac{\partial^2}{\partial s \partial t}}
P(s,t;u)\right|_{(0,0)} = -\int_{0}^{u} \left( a(\nu)
\left. {\tfrac{\partial}{\partial t}} P(0,t;\nu) \right|_{t=0} +
 b(\nu) 
\left.{\tfrac{\partial }{ \partial s}} P(s,0;\nu)\right|_{s=0}\right)
d\nu.$$
The equations
$\left. \frac{\partial}{\partial t} P(0,t;\nu) \right|_{t=0}
= - \int_0^\nu b(\mu) d \mu$
and
$\left.\frac{\partial }{ \partial s}P(s,0;\nu)\right|_{s=0}
= - \int_0^\nu a(\mu) d \mu$
can be obtained from (\ref{curly}) in a similar manner,
using that $P(0,0;u)$ is the identity.
Substituting each of these into the equation above and evaluating at $u=1$
gives the desired result since
$\hol(A+sa+tb)=\hol(A)P(s,t;1)$.
\end{pf}

This proposition allows us to compute
the first and second derivatives of any function $f:\AA \lto \RR$
of the form $f= \tau \circ \hol,$
where $\tau:SU(3) \lto \RR$ is a smooth invariant function.
An important example is when $\tau$ is either
the real or imaginary part of $\tr:SU(3) \lto \CC.$

\begin{cor}\label{derivative of tr hol}
The first and second derivatives of the trace of holonomy are given by:
\begin{itemize}
\item [(i)] $\left. \frac{d}{dt }
\tr(\hol(A+ta))\right|_{t=0}=
\int _{0}^{1} \tr \left( \hol (A)  a(\mu)  \right)d\mu,$
\item [(ii)] $\left. \frac{\partial^{2}}{\partial s \partial t}
\tr(\hol(A+sa+tb))\right|_{(0,0)}=
\int _{0}^{1}\int _{0}^{\nu} \tr \{\hol (A)  (a(\nu)  b(\mu) +b(\nu)a(\mu))
\}d\mu
d\nu.$
\end{itemize}
\end{cor}
\noindent{\it Remark.}  Proposition \ref{derivative of holonomy} and Corollary 
\ref{derivative of tr hol} remain valid for $SU(n), \, n>3.$ 

In Section 3, we shall  show that for a suitable choice of $\Ga$, 
regularity of $\MM_{h}$ is a generic condition for 
$h\in \ps_{\Ga}$ near zero, and similarly for 
regularity of $W_{\rho}$ for $\rho \in C^{1}([-1,1], 
\ps_{\Ga})$.
The following proposition provides useful bounds on the derivatives
of admissible functions.

\begin{prop}\label{bounds from taubes}
\begin{enumerate}\item[(i)] Fix  $\ga:S^{1}\times D^{2} \lto X$ an embedding
of the solid torus and 
let $\tau_1,\tau_2$ be the real and imaginary parts of trace on $SU(3).$
Then there exists a constant $C_1$ depending on $\ga$ 
such that
$$ |D^{n}p(\ga, \tau_j) (A)(a_{1},\ldots, a_{n})|\leq C_{1} 
\prod_{i=1}^{n} \| a_{i}\|_{L^{2}_{1}}$$ 
for all $A \in \AA$ and for $j=1,2.$ 
\item[(ii)]
Fix $\Ga$ a collection of
embedded solid tori. Then there exists a constant $C_{2}$
depending on $\Ga$ such that
the inequalities hold for all $h\in \ps_{\Ga}$ and all $A \in \AA$
\begin{eqnarray*}
|Dh(A)(a_{1})| &\leq&  C_{2} \, \|h\|_{C^3} \cdot \| a_{1}\|_{L^{2}_{1}}, \\
|D^2 h(A)(a_{1},a_2)| &\leq&  C_{2} \,
\|h\|_{C^3} \cdot \| a_{1}\|_{L^{2}} \cdot \| a_2 \|_{L^2}, \\
|D^3 h(A)(a_{1},a_2,a_3)| &\leq&  C_{2} \, 
\|h\|_{C^3} \cdot \| a_{1}\|_{L^{2}_{1}} \cdot \| a_2 \|_{L^2}\cdot  \| a_3\|_{L^2}, \\
\| \nabla h(A) \|_{L^{2}_{1}} &\leq& C_{2}\,  \|h\|_{C^3}. 
\end{eqnarray*}
\end{enumerate}
\end{prop}

\begin{pf}
See \cite{taubes}, Section 8a.\end{pf}

The last  proposition of this section
allows  one to patch together
the local regularity arguments to give global results in subsection \ref{abs}.
\begin{prop} \label{compactness}
If $C\subset \ps$ is compact, then
$\bigcup_{h\in C} \MM_{h}$
is also compact.
\end{prop}\begin{pf}  See Lemma 8.3 in \cite{taubes}.\end{pf}

\newpage
\section{Transversality}
\label{transversality section}
The goal of this section is to  establish various structure theorems for
the perturbed flat moduli space $\MM_h$
and for the parameterized moduli space $W_\rho$
for generic  $h \in \ps$ and generic  $\rho \in C^1([-1,1],\ps).$
Before doing this, we must fix a collection
$\Ga$ of solid tori so that 
the resulting 
space of perturbations $\ps_\Ga$ 
is general enough for these transversality results to hold.

The first subsection contains a
formulation 
of the necessary conditions on $\Ga$ 
and a result which implies that we can always choose
$\Ga$ to satisfy these conditions in a neighborhood of $\MM$
in $\BB \times \ps_\Ga.$
In the second subsection, we proceed with
the transversality results for $\MM_h$ and $W_\rho$.

\subsection{Abundance of admissible functions}
\label{abs}

For any $A \in \AA,$ define 
$$\KK_A = \ker d_A^* \cap \Om^1(X;su(3))$$ and denote by
$\Pi_A: \Om^1(X;su(3)) \lto \KK_A$
the $L^2$ orthogonal projection.
The slice through $A$ to the gauge action is the affine subspace
$$ X_{A}=\{ A+ a \mid a\in \KK_A \} \subset \AA.$$
A small neighborhood of $A$ in $X_{A}$, divided by the stabilizer of
$A$,  gives a local model for
$\BB$ near $[A]$.

The first proposition reduces the study of the local structure of the
moduli space to a Fredholm problem.
\begin{prop}\label{make fredholm}
Given a perturbed flat connection, there is a neighborhood $U\subset
X_{A}$ of $A$ such that $A+a\in U$ implies that
$\zeta_{h}(A+a)=0$ if and only if
$\Pi_{A}\zeta_{h}(A+a)$.
\end{prop}
\begin{pf}  See Lemma 12.1.2 of \cite{mmr} and Lemmas 28 and 29 of
\cite{herald1}.\end{pf}

\begin{defn} \label{hermite}
Suppose $A$ is a  reducible h-perturbed flat connection and denote by
$\sym \HH^{1}_{A,h}(X;\rrp)$ the set of $\stab (A) \cong U(1)$ invariant 
symmetric (hence Hermitian) bilinear forms on
$\HH^{1 }_{A,h}(X;\rrp)$. 
\end{defn}  

\begin{defn}  \label{abundance def}
A collection $\Ga$ of embedded solid tori in $X$
is called {\bf abundant}  for
$(A,h),$ where
$h\in \ps_{\Ga}$ and $A \in \AA^* \cup \AA^r$ is
h-perturbed flat, in case
there exists a finite subset $\{f_1,\ldots, f_m\} 
\subset \ps_\Ga$ of admissible functions such that:
\begin{enumerate}
\item[(i)] If $ A \in \AA^*,$ then 
the map from $ \RR^{m} $ to $\hom(\HH^1_{A,h}(X;su(3)), \RR)$ given by \\
$(x_1,\ldots, x_{m}) \mapsto \sum_{i=1}^{m} x_i Df_{i}(A)$ is surjective.
\item[(ii)] If $A \in \AA^r,$ then 
the map from
$\RR^m$ to $\hom(\HH^{1}_{A,h}(X; \rr),\RR)\oplus \sym
\HH^{1}_{A,h}(X;\rr^\perp)$ given by
$(x_1,\ldots, x_m) \mapsto \left( \sum_{i=1}^{m}x_{i}Df_{i}(A), \sum_{i=1}
^{m}
x_i \hess f_i(A) \right)  $
is surjective.  
\end{enumerate}
Because abundance is a gauge invariant concept, it makes sense
to say that $\Ga$ is abundant for $([A],h)$. 
When $h=0,$ we say that
$\Ga$ is abundant for $A$ or $[A]$.
\end{defn}

If $\Ga$ is abundant for $(A,h)$ and $\Ga \subset \Ga',$
then of course $\Ga'$ is also abundant for $(A,h).$
The next proposition is the principal result of this subsection;
it shows that there exists a collection $\Ga$ which is
abundant for all nontrivial perturbed flat connections in a neighborhood
of the flat moduli space. This is a global result
and its proof will occupy the remainder of the subsection.
The statement of the proposition is divided into three parts, 
which can be viewed as the
pointwise, local, and global versions of the same result.

\begin{prop}\label{openness}
\begin{enumerate}
\item[(i)] If $A \in \AA$ is a nontrivial flat connection,
then there exists a finite collection
$\Ga$ which is abundant for $A.$
In case $A$ is reducible, $\Gamma$ and the 
subset $\{f_1, \dots, f_m\}$ from Definition \ref{abundance def}
can be chosen so that for some $k,$ 
\begin{enumerate}
\item[(a)] $\{ Df_1(A),\ldots, Df_k(A)\}$ spans 
$\hom(\HH^{1}_{A}(X; \rr),\RR)$
\item[(b)] $\{ \hess f_{k+1}(A), \ldots, \hess f_m(A) \}$ spans 
$\sym \HH^{1}_{A}(X;\rr^\perp)$
\item[(c)] $Df_j(A) = 0$ for $j=k+1,\ldots, m.$
\end{enumerate}

\item[(ii)] If $A \in \AA$ is a nontrivial flat connection
and $\Ga$ is abundant for $A$ and is chosen as in (i),
then there exists an
open neighborhood $U\times V\subset \BB \times \ps_\Ga$ of $([A],0)$ 
such that $\Ga$ is abundant 
for all $([A'],h)\in U\times V$ with $\zeta_{h}(A')=0$.
\item[(iii)]
There exist a finite collection $\Ga$ and an open
neighborhood $U\times V\subset \BB \times \ps_\Ga$ of
$\MM \setminus [\th]$ such that $\Ga$ is abundant 
for all $([A],h)\in U\times V$ with $\zeta_{h}(A)=0$.
\end{enumerate}
\end{prop}
\begin{pf} Part (i) follows from Corollary \ref{florida} and
Proposition \ref{hessian loops}, as we now explain.
First, assume $A$ is irreducible.  Replace all loops
$\ell$ coming from \ref{florida} (ii) 
by tubular neighborhoods $\ga.$
Next, by shrinking the tubular neighborhoods, if
necessary, we can approximate
functions $f:\AA \lto \CC$ of the form $f(A) = \tr(\hol_{\ell}(A))$ 
arbitrarily closely 
by the complex-valued functions $p(\ga,\tr)(A)$ defined 
as in equation (\ref{elementary function}). 
In case $A$ is reducible, apply the same procedure to obtain
real-valued functions $p(\ga, \tr_\RR)(A)$ from the
real part of $\tr(\hol_{\ell}(A))$ for the loops
in \ref{florida} (i). 
This proves (i) for $A$ irreducible as well as part (a)
for $A$ reducible.

To finish off part (i) in case $A$ is reducible,
thicken the loops obtained from an application of Proposition
\ref{hessian loops}. This provides a collection of functions
with $Dp(\gamma,\tr)(A)=0$ 
whose Hessians span
$\sym \HH^{1}_{A,h}(X;\rr^\perp)$. This proves (b) and (c) and completes
the proof of part (i).

Part (ii) says that abundance is an open condition 
around flat connections in $\AA \times \ps_\Ga$
and requires several estimates,
contained in Lemmas \ref{first lem} and \ref{las vegas}. 
Before presenting those arguments,
we explain how (iii) follows from (i) and (ii).
  
By (i) and (ii), for any nontrivial flat connection $A,$
we have a collection $\Ga$ which 
is abundant for all perturbed flat orbits $([A'],h)$ 
in a neighborhood $U' \times V' \subset \BB \times \ps_{\Ga}$ 
of $([A],0)$. Applying this for each $[A] \in \MM \setminus [\th]$ and 
using compactness, we obtain a finite subcover
$U'_1, \ldots, U'_{l}$ 
and corresponding collections $\Ga_1, \ldots, \Ga_l$. 
Set $\Ga = \bigcup_{i=1}^{l} \Ga_i$.  
Part (iii) follows by applying (ii) once again to
$A$ and the collection
$\Ga$ 
to obtain an open neighborhood $U \times V \subset \BB \times \ps_{\Ga}$ 
of $([A],0)$  
such that $\Ga$ is abundant for
all $([A'],h) \in U \times V$ with $\zeta_h(A')=0.$
This last step is performed for each $[A] \in \MM \setminus [\th]$,
and compactness once again allows us to extract a finite subcover
$U_1, \ldots, U_k$ of $\MM \setminus [\th]$.
The proof of part (iii) is completed by
setting $U = \bigcup_{i=1}^k U_i$
and $V = \bigcap_{i=1}^k V_i.$

As for part (ii), it is easiest to
see this in case $A$ is irreducible.
On the other hand, if $A$ is reducible,
then similar reasoning shows that abundance is local
in $\BB^{r}\times \ps_\Ga,$  
but whether there exists an open neighborhood in 
$\BB \times \ps_\Ga$ is less obvious.
The following argument treats
{\it irreducible} perturbed flat connections 
in a neighborhood of $A$ assuming $A$ is reducible.
Before continuing with the proof,
we need to introduce some notation.

Since $A$ is a fixed
reducible flat connection
for the rest of this proof,  
we write $\KK$ for $\KK_A.$
It is useful to decompose elements $ a \in \KK$ 
as $a = (a_1,a_2)$ according to 
$su(3) = \rr \oplus \rrp$. 
Thus $a_1 \in \Om^1(X;\rr)$ and $a_2 \in \Om^1(X;\rrp).$
For $i=1,2$, we have the Hodge decomposition $a_i = (a_i', a_i'')$ 
where $a_1' \in \HH^1_A(X;\rr)$ and $a_2' \in \HH^1_A(X;\rrp)$
are the cohomological components
and $a_1'', a''_2$  are characterized as follows.
Define $\KK''_{1}$ 
to be the orthogonal complement of
$\HH^1_A(X;\rr)$ in $\KK \cap \Om^1(X;\rr)$,
and also $\KK''_{2}$ to be the orthogonal complement of
$\HH^1_A(X;\rrp)$ 
in $\KK \cap \Om^1(X;\rrp)$.
Denote by $\Pi_i'': \Om^1(X;su(3)) \lto \KK_i''$
the $L^2$ orthogonal projection 
for $i=1,2.$
Then 
$a_i'' = \Pi''_i a  \in \KK''_i$ 
and $a = (a_1,a_2) = (a_1',a_1'',a_2', a_2'').$
We set $\KK'' = \KK''_1 \oplus \KK''_2$ 
and $\Pi'' = (\Pi_1'', \Pi_2'').$ 

Suppose $a,b \in \Om^1(X;su(3)).$ The notation $[a \wedge b]$ 
indicates the
product obtained by combining the wedge product on the form
part with the Lie bracket on the coefficients.
The following is the $su(3)$
analog of the
well-known formulas for the Lie bracket in $su(2)$ 
(with regard to the decomposition $su(2) = u(1) \oplus u(1)^\perp$).
If we decompose $a=(a_1,a_2)$ 
and $b=(b_1,b_2)$ according to $su(3) = \rr \oplus \rrp$ as above,
then 
$$\left\{ \begin{array}{cc}
*[a_i \wedge b_j] \in \Om^1(X;\rr) & \hbox{ if } i=j,\\
*[a_i \wedge b_j] \in \Om^1(X;\rrp) & \hbox{ if } i \neq j. 
\end{array} \right.$$

The proof proceeds with
two lemmas.  
The  first one shows
that the space of perturbed flat irreducible connections in $X_{A}$
for small $h$ are close to the image of the affine subspace
$A + \KK''_1 + \HH^1_A(X;su(3)).$
It also gives some control over the
distance from the nearby reducibles to the affine subspace
$A+\HH^{1}_{A}(X;su(3))$ in terms of the size of the perturbation.

\begin{lem}\label{first lem}  For any $\Ga$ and any $0<R<1$,
there exist
$K<\infty$ and  $0<\ep<1$ such that  if $A+a\in X_{A}$  is $h$-perturbed flat
with $\|a\|_{L^{2}_{1}}<\ep$ and $\|h\|_{C^3}<\ep$, then
\begin{itemize}
\item [(i)] $\| a''_{2}\|_{L^{2}_{1}} \leq R \, \|a'_{2}\|_{L^{2}_{1}}$  and
\item [(ii)] $\|
a''_{1} \|_{L^{2}_{1}} \leq R \left( \| a'_{1}\|_{L^{2}_{1}} + \|a'_{2} \|_{L^{2}_{1}}\right) + K\, \| h \|_{C^3}$.
\end{itemize}
\end{lem}
\begin{pf}
Fix $0<R<1$.
Consider the map from $X_{A}\times \ps$ to
$\KK''$ given by $\Pi'' \zeta_{h}(A+a).$ 
The linearization at $(A,0)$  restricted to $\KK''$
with the $L^{2}_{1}$  
norm on the domain and $L^{2}$ norm on the range
is $*d_{A}$, an elliptic Fredholm operator with 
trivial kernel.
Therefore there exists $\la>0$ such that 
$\|*d_{A}b''\|_{L^{2}}\geq \la \| b'' \|_{L^{2}_{1}}$ 
for all $b''\in \KK''$.

Now assume that $\Pi_{A}\zeta_{h}(A +a)=0$.  Expanding the equation
$\Pi''_{2} \, \zeta_{h}(A+a) = 0$ gives
$$
0=*d_{A}(a''_{2}) + 2 \Pi_{A} *[a_{1}\wedge a_{2}] 
- 4 \pi^2 \Pi''_2 \, \nabla h(A+a).$$
By Taylor's theorem, the last term on the right can be replaced by 
$$- 4 \pi^2 \left[ \Pi''_2 \left( \hess h(A+a_{1})(a_{2}) 
+  D^{2}\nabla h (A+a_{1} + t_1 a_{2})(a_{2}, a_{2}) \right) \right],$$
for some $0<t_1<1$.
Here we are exploiting the equivariance of $\zeta_h$ with respect to the
$\stab(A) $ action.
Rearranging and using the triangle inequality
on $a_2 = a_2'+a_2''$, we obtain
$$\la \| a''_{2}\|_{L^{2}_{1}} \leq \left( 2C \|  a_{1}\|_{L^{2}_{1}}  
+ 8 \pi^2 C_{2}\| h \|_{C^3} \right) \left(\| a'_{2}\|_{L^{2}_{1}} +
\| a''_{2}\|_{L^{2}_{1}}\right),$$
where $C$ comes from the Sobolev multiplication theorems and $C_{2}$
is the constant given in Proposition \ref{bounds from taubes}.
By shrinking $\ep$ to control some of the $L^{2}_{1}$ norms on the
right side, we obtain the first claim.

To prove the second claim, expand the equation $0=\Pi''_1 \, \zeta_{h}(A+a)$ to get
$$
0=*d_{A}(a''_{1}) 
+ \Pi_{A} * \left( [a_{1}\wedge a_{1}]+ [a_2 \wedge a_{2}]\right)
 - 4 \pi^2 \Pi''_1 \, \nabla h(A+a).$$
Rearranging, we see that
$$\la \| a''_{1}\|_{L^{2}_{1}} \leq 
C \left(\| a_{1} \|_{L^2_{1}}^{2} + \| a_{2}\|_{L^{2}_{1}}^2 \right) 
+ 4 \pi^2 C_{2}\| h \|_{C^3}.$$
Now apply the triangle inequality on the right to $a_1 = a_1' + a_1''$
and use the first part to obtain the required bound.
\end{pf}

The next lemma is a similar result about tangent vectors at
perturbed flat connections which are in the kernel of the Hessian
of $\CS + h$ (restricted to $X_{A}$).  We decompose $b\in T_{A+a}X_{A}$ into
$b=( b_1,b_2) = (b_1', b_1'', b_2', b_2'')$
as before.

\begin{lem}  \label{las vegas}  For any $\Ga$ and any $0<R<1$, there
exist $K<\infty $ and $0<\ep <
1$ such that
if $A+a\in X_{A}$ is a nonabelian $h$-perturbed flat
with $\|a\|_{L^{2}_{1}}<\ep$ and $\| h \|_{C^3}<\ep$, and if $b\in T_{A+a}X_{A}$
is in the kernel of $\hess (\CS + h)(A+a)$, then
\begin{itemize}
\item [(i)] $\| b''_{1}\|_{L^{2}_{1}} < R\| b'_{1} \|_{L^{2}_{1}} +
K \| a'_{2}\| _{L^{2}_{1}} \cdot \| b'_{2}\| _{L^{2}_{1}}$
\item [(ii)] $\| b''_{2}\|_{L^{2}_{1}} < R\| b'_{2} \|_{L^{2}_{1}} +
K \| a'_{2}\| _{L^{2}_{1}} \cdot \| b'_{1}\| _{L^{2}_{1}}$
\end{itemize}
\end{lem}
\begin{pf}
Setting the $\rr$ and $\rrp$ components of
$D \, \Pi_{A}\zeta_{h}(A+a)(b)$ equal to zero gives two coupled equations in
$b_{1}$ and $b_{2}$.  Expanding the $\rr$
component leads to
\begin{eqnarray*}
-*d_{A}b''_{1}&=& \Pi_{A} *\left( [a_1 \wedge b_{1}] + [a_{2}\wedge b_{2}] \right)
-4 \pi^2   \Pi_{A}\hess h(A+a_{1})(b_{1}) \\
&&- 4 \pi^2  \Pi''_1 D \{\hess h(A+a_{1}+t_{2}a_{2})(b)\} (a_{2}).
\end{eqnarray*}
Taking the $L^{2}$ norm of each side of this equation 
and using the various bounds
as in the last lemma, it follows that
$$\la \| b''_{1} \|_{L^{2}_{1}} \leq
C \left(\| a_{1} \|_{L^{2}_{1}} \cdot \| b_{1}\| _{L^{2}_{1}}
+\|a_{2}\|_{L^{2}_{1}} \cdot \| b_{2}\| _{L^{2}_{1}}\right)
+ K \, \| h \|_{C^3} \left(\| b_{1} \|_{L^{2}_{1}}
+\| a_{2}\|_{L^{2}_{1}}
\cdot \| b \|_{L^{2}_{1}}\right).$$
Applying the triangle inequality, first to $b=b_1+b_2$ and
then  to $ b_1 = b_1'+b_1''$ 
everywhere on the right hand side of this 
equation and moving all occurrences of $b''_{1}$ to the left,
we see that, for $\ep$ small enough, 
\begin{eqnarray}
\frac{ \la}{2} \| b''_{1}\|_{L^{2}_{1}} & \leq &
2 C \left( \ep \, \| b'_{1} \|_{L^{2}_{1}}  + \| a'_{2} \|_{L^{2}_{1}}
\cdot \| b_{2} \|_{L^{2}_{1}} \right)
+  K \ep \left( \| b'_{1}\|_{L^{2}_{1}}+
2 \| a'_{2}\|_{L^{2}_{1}} \cdot
\| b_{1}' +  b_{2} \|_{L^{2}_{1}} \right)  \nonumber \\
&\leq&  
\ep \, \const \| b'_{1} \|_{L^{2}_{1}}  
+ \const \| a'_{2} \|_{L^{2}_{1}} \cdot
\| b_{2} \|_{L^{2}_{1}}. \label{fourtharray}
\end{eqnarray}

Similar reasoning applied to the $\rrp$ component of
$D \, \Pi_A \zeta_h(A+a)(b)$ gives
\begin{equation}
\label{fiftharray}
\frac{\la}{2} \, \| b''_{2}\|_{L^{2}_{1}} \leq  \ep \, \const 
\| b'_{2} \|_{L^{2}_{1}}  + \const \| a'_{2} \|_{L^{2}_{1}}\cdot 
\|b_{1} \|_{L^{2}_{1}} .
\end{equation}
The conclusion of the lemma follows
from equations (\ref{fourtharray}) and (\ref{fiftharray}).
\end{pf}

We are now ready to complete the proof of Proposition \ref{openness} (ii).
Referring to part (i), since $A$ is reducible,
we have finite subsets $\{ f_1, \ldots, f_k \}$ and $\{ g_1, \ldots,
g_l\}$ 
of $\ps_\Ga$ such that  
\begin{itemize}
\item [(i)] $\Span \{ Df_{i}|_{\HH^{1}_{A}(X;\rr)} \mid i=1,\ldots, k\}
= \mbox{Hom} (\HH^{1}_{A}(X;\rr),\RR),$
\item [(ii)] $
\Span \{ D^{2} g_{j}|_{\HH^{1}_{A}(X;\rrp)^{\otimes 2}} \mid j=1,\ldots, l\} 
=\sym \HH^{1}_{A}(X;\rrp )$.
\item [(iii)] $D g_j|_{\HH^{1}_{A}(X;\rr)} =0$ for $j=1, \ldots, l$.
\end{itemize}
Our strategy here is to show that, given $a$ and $h$
sufficiently small with $A+a$ an irreducible $h$-perturbed flat connection,
the functions $\{f_i, g_j \}$ 
detect all elements $b \in \ker K(A+a,h)$  to first order.

Choose a constant $N>0$  such that, for 
all $u\in \HH^{1}_{A}(X;\rr)$ and all $v,w \in
\HH^{1}_{A}(X;\rrp)$, the following
bounds hold:
\begin{eqnarray}
\max_{1\leq i \leq k} \{ | Df_{i} (A) (u) | \} 
&\geq & N  \, \| u \|_{L^{2}_{1}}\label{f N} \\
\max_{1\leq j \leq l}  \{ | D^{2}g_{j}(A)(v,w) |\} 
&\geq & N  \, \|v\|_{L^{2}_{1}}  \cdot \|w\|_{L^{2}_{1}}\label{g N}.
\end{eqnarray}
Choose $\ep $ small enough that these inequalities  
continue to hold when $N$ is replaced by $\frac{N}{2}$
and $A$ is replaced by $A+a$ for $\| a \|_{L^{2}_{1}} < \ep$.

Suppose that $h\in \ps_{\Ga}$ and that 
$A+a\in X_{A}$ is an irreducible h-perturbed flat connection,
and assume $b \in \Om^1(X;su(3))$
is an element in the kernel of $\hess (\CS + h)(A+a)$.
Choose  functions $f$ and $g$ from $\{f_{i} \}$ and $\{g_{j}\}$,
respectively, for which
$|Df(A+a)(b'_{1}) |\geq N/2 \, \| b'_{1}\|_{L^{2}_{1}}$ and
$|D^{2}g(A+a)( a'_{2}, b'_{2})|\geq N/2 \, \| a'_{2}\|_{L^{2}_{1}} \cdot
\| b'_{2} \|_{L^{2}_{1}}$.
If either
$Df(A+a)(b)$ and $Dg(A+a)(b)$ is non-zero,
then we are done. So we assume both vanish 
and seek a contradiction.

Apply the triangle inequality to the equation
$Df(A+a)(b'_{1})=-Df(A+a)(b''_{1}+b_{2})$ to get the
inequality
\begin{eqnarray*}
\frac N 2 \| b'_{1} \|_{L^{2}_{1}}& \leq &| Df(A+a)(b''_{1}) |
	+ |Df(A+a_{1})(b_{2})|  + |D^{2}f(A+a_{1}) (a_{2}, b_{2})| \\
&+& |D^{3}f(A+a_{1}+ t_{1}a_{2})(a_{2}, a_{2}, b_{2})|, 
\end{eqnarray*}
where $0<t_{1}<1$.
Then  $Df(A+a_{1})(b_{2})$ is zero by invariance
under $\stab (A+a_1) \cong U(1)$,   
and applying bounds to the other
terms gives
\begin{equation}
\frac N 2 \| b'_{1} \|_{L^{2}_{1}} \leq C_{2} \, \| f \|_{C^3} \cdot \| b''_{1}\|_{L^{2}_{1}}
+ 4C_{2} \, \| f \|_{C^3} \cdot \| a'_{2} \|_{L^{2}_{1}} \cdot \| b_{2} \|_{L^{2}_{1}}
\nonumber 
\end{equation}

Using Lemma \ref{las vegas}, and choosing  $\ep$ suitably
small, this implies
\begin{equation}
\label{third f bound}  
\frac{N}{3}  \, \| b'_{1} \|_{L^{2}_{1}} \leq
\const \, \| a'_{2}\|_{L^{2}_{1}} \cdot  \| b'_{2} \|_{L^{2}_{1}}.
\end{equation}

Next consider $Dg(A+a)(b)$.  We first bound the derivative in the
$b_{1}$ direction. 
\begin{eqnarray}
|Dg(A+a)(b_1)| &=& |Dg(A)(b_{1}) +D^{2}g(A+t_{1}a)(a_{1}, b_{1})
+D^{2}g (A+t_{2}a)(a_{2}, b_{1}) | \nonumber \\
&=&|D^{2}g(A+t_{1}a)(a_{1}, b_{1})
+ D^{2}g(A+t_{2}a_{1})(a_{2}, b_{1})
+ D^{3}g(A_{1}) (t_2 a_{2}, a_{2}, b_{1})| 
\nonumber \\
&\leq& C_{2} \| g \|_{C^3}\cdot \|b_{1}\|_{L^{2}_{1}} 
\left( \| a_{1} \|_{L^{2}_{1}}
+ \|  a_{2} \|^2_{L^{2}_{1}} \right)  
\leq \ep \, C_3 \| b_{1}\|_{L^{2}_{1}} 
\nonumber \\
&\leq & \ep \, \const  \, \| b_1' \|_{L^{2}_{1}} + \ep \, \const \, 
\| a'_{2}\|_{L^{2}_{1}} \cdot \| b'_{2} \|_{L^{2}_{1}}.
\label{first g bound}
\end{eqnarray}
In the first line, 
$Dg(A)(b_{1})=0$ by hypothesis,
and in the second, $D^{2}g(A+t_{2}a_{1})(a_2,b_1)$
vanishes by gauge symmetry.
The last step follows from part (i) of Lemma \ref{las vegas}.

Finally, we bound the derivative of $g$ in the $b_{2}$ direction away
from zero.
$$|Dg(A+a)(b_{2})|=|Dg(A+a_{1})(b_2) + D^{2}g(A+a_{1})(a_{2}, b_{2}) 
+ D^{3}g(A_{2})(a_{2}, a_{2}, b_{2})| $$
Appling gauge symmetry once more shows that
$Dg(A+a_{1})(b_2)=0$ in the equation above.
Bounds on the other terms give, for $\ep $
sufficiently small,
\begin{eqnarray}
|Dg(A+a)(b_{2})| & \geq & \frac{N}{2} \|a'_{2}\|_{L^{2}_{1}} \cdot
\|b'_{2}\|_{L^{2}_{1}} - \const \| a''_{2}\|_{L^{2}_{1}} \cdot
\|b'_{2} \|_{L^{2}_{1}} \nonumber \\
&& - \; \const \|a''_{2}\|_{L^{2}_{1}} \cdot \| b''_{2}\|_{L^{2}_{1}} 
- \const \| a'_{2} \|_{L^{2}_{1}} \cdot \| b''_{2} \|_{L^{2}_{1}} \nonumber \\
& \geq & \frac{N}{3}\| a'_{2} \|_{L^{2}_{1}} \cdot \| b'_{2} \|_{L^{2}_{1}}
- \const  \| a'_{2}\|^{2}_{L^{2}_{1}} \cdot \| b'_{1}\|_{L^{2}_{1}}
\label{second g bound}
\end{eqnarray}

Combining inequalities (\ref{first g bound}) and 
(\ref{second g bound}), we get
$$\frac{N}{4}\| a'_{2}\|_{L^{2}_{1}} \cdot \| b'_{2}\|_{L^{2}_{1}}
\leq \ep \, \const \| b'_{1} \|_{L^{2}_{1}},$$
which, combined with inequality (\ref{third f bound}), gives the
desired contradiction. 
\end{pf}

Since $X$ is an integral homology 3-sphere, there are no noncentral
abelian flat connections.  The following proposition guarantees that
this, together with the property that
${\Ga}$  is abundant, continue to hold for small perturbations.
It also provides a unique component of the flat moduli space near each
perturbed flat connection, for small perturbations.

\begin{prop}\label{existence of epsilon}
Suppose $\Ga$ satisfies condition (iii) of Proposition
\ref{openness}. There exists an $\ep_{0}>0$
such that:
\begin{enumerate}
\item[(i)] If $A \in \AA^*\cup \AA^{r}$ is flat and $A' \in \AA$ is abelian,
then $\|A-A'\|_{L^{2}_{1}}>2\ep_{0}.$
\item[(ii)]  If $\|h\|_{C^{3}}< \ep_{0}$ and $A \in \AA$ is
$h$-perturbed flat, then there exists $\widehat{A} \in \AA$ which is flat
with $\|A- \widehat{A}\|_{L^{2}_{1}} < \ep_{0}.$
\item [(iii)]  If $\|h\|_{C^{3}}< \ep_{0}$ and $A \in \AA$ is
$h$-perturbed flat, then ${\Ga}$ is abundant for $([A],h)$.
 \item[(iv)] If $A,A' \in \AA^r$ are flat and lie on
different components of the space of flat connections in $\AA$,
then $\| A - A'\|_{L^{2}_{1}} > 2 \ep_{0}.$
\end{enumerate}
\end{prop}
\begin{pf}  For claims (i) and (ii), see
Lemma 1.3 and Proposition 1.5 of \cite{taubes}.
Claim (iii) follows from claim (ii).  For the neighborhoods $U$ and
$V$ in Proposition \ref{openness},
choose $\ep_{0}$ small enough that
the ball of  radius $\ep_{0}$ around $0\in \ps_{\Ga}$
is contained in  $V$
and the $\ep_{0}$ neighborhood of $\MM^{*}\cup \MM^{r}$ in $\BB$
is contained in $U$.

For part (iv), suppose to the contrary that 
there were no $\ep_0$ 
satisfying the conclusion.
Then we have two sequences 
$A_i$ and $A' _i$ of flat connections in $\AA$
with
$\| A_i -A'_i\|_{L^2 _1	}<\frac 1 i$ such that
$A_i$ and $A'_i$ never lie on
the same component of the space of flat connections.
By compactness of $\MM,$
after passing to a subsequence, we can assume that
there is a sequence of gauge transformations $g_i$ such that
$g_i \cdot A_i$ converges to a flat connection
$A_0$.  Then $g_i \cdot A'_i$ must also converge to $A_0$. (Note that
we are using the
standard gauge invariant $L^2 _1$ norm here.)

Consequently, 
for $i$ large, we see that
$g_i \cdot A_i$ and $g_i \cdot A'_i$
must lie on
the same component of the space of flat connections
as the one containing $A_0.$
But this implies that $A_i$
and $A'_i$ lie on the
same component, which is a contradiction.
\end{pf}

\subsection{Regularity theorems}
\label{regularity theorems}
We are now ready to prove the structure theorems for $\MM_h$
and $W_\rho.$
We begin with the definition of regularity in this context.
Throughout this subsection,
$\Ga$ denotes a fixed collection of solid tori satisfying
Proposition \ref{openness}, part (iii).
Thus, $\Ga$ is abundant for all pairs
$([A],h) \in \BB \times \ps_\Ga$ in a neighborhood of 
$(\MM \setminus [\th]) \times \{0\}$.
\renewcommand{\ps}{{{\cal F}(\ep_0)}}
Choose $\ep_{0}$ as in Proposition
\ref{existence of epsilon} and define $\ps$ to be
the ball of radius $\ep_{0}$ about $0$ in the space ${\cal F}_\Ga$
of admissible functions.

\begin{defn} \label{reg} Suppose $h\in \ps$ and
 $U\subset\MM_{h}$ is open. Then $U$ is {\bf regular} in case
$\HH^{1}_{A,h}(X;su(3))$ is trivial for all $[A]\in U$.
\end{defn}
Regularity as defined here makes no assumption
on the irreducibility of $A.$
\begin{prop}  \label{regularity1} If $U\subset\MM _{h}$
is regular, then $\MM_{h} ^{*}\cap U$ and $\MM_{h}^{r}\cap U$ are
0-dimensional submanifolds of $\BB^{*}$ and $\BB^{r}$,
respectively.
\end{prop}
\begin{pf}
This follows directly from standard Kuranishi arguments.  \end{pf}

We define regularity for the parameterized moduli space next. 
For any triple $(A,\rho, t)\in \AA\times C^{1}([-1,1],\ps) \times
[-1,1]$,
define an index one Fredholm operator
by the formula
\begin{eqnarray*}
&L(A,\rho,t)\, :\, \Om^{0+1}(X;su(3))\oplus \RR \lto
\Om^{0+1}(X;su(3))& \\
&(\xi , a, \tau) \mapsto K(A,\rho_t)(\xi, a) - 4 \pi^2 \tau
{\tfrac{\partial}{\partial t}} \nabla \rho_t (A)
\end{eqnarray*}

Since $X$ is an
integral homology 3-sphere, the only abelian
orbit in the flat moduli space
is $[\th],$ and this continues to be true 
for small perturbations
thanks to Proposition \ref{existence of epsilon}.
This explains why we dismiss the case of abelian orbits
in the following definition.
Note, however, that such orbits may indeed occur
for large perturbations, or even for small perturbations
on arbitrary 3-manifolds.

\begin{defn} \label{p m reg}   Let $\rho
:[-1,1]\lto \ps$ be a $C^{1}$ curve with $\MM_{\rho(\pm 1)}$
regular. An open subset $U\subset W_{\rho}$ is
{\bf regular} if:
\begin{enumerate}
\item[(i)] $\HH^{1}_{\th, \rho_t}(X;su(3))$ is trivial for
$([\th], t) \in U$.
\item[(ii)]  $U$ contains no noncentral abelian orbits.
\item[(iii)] For all $([A], t)\in W_{\rho}^{*}\cap U$, $ L(A,\rho, t)$ is
surjective.
\item[(iv)] For all $([A], t)\in W_{\rho}^{r}\cap U$, 
$ \Om^{0+1}(X; \rr) \cap \coker L(A,\rho,t)=\HH^{0}_{A}(X;su(3))\cong u(1)$.
\item[(v)] There is a  finite
subset $J$ of $W_{\rho}^{r}\cap U$ such that
for $([A],t) \in W^r_\rho,$
$$\dim \HH^{1}_{A,\rho_t}(X;\rrp)= \left\{ \begin{array}{ll} 2 & \hbox{if }
([A],t) \in J \\
   0 & \hbox{otherwise.} \end{array} \right.$$
Elements of $J$ are called {\bf bifurcation points}.
 \item[(vi)] If $([A_s],t_s)$
is a parameterized curve in  $W_{\rho}^{r} \cap U$  and $([A_0],t_0)
\in J$, then the
(multiplicity two) eigenvalue of $K(A_s,\rho(t_s))$ crosses
zero transversally at $s=0.$
\end{enumerate}\end{defn}
Note that regularity of $W_\rho$ does not ensure that
$W_\rho$ is a smooth cobordism (cf.  Lemma \ref{regular cobordism}).
Conditions (v) and (vi) of Definition \ref{p m reg} make sense in light of
claim (i) of the next lemma.
\begin{lem}  \label{regular cobordism}
If $U^{r}\subset W_{\rho}^{r}$ is open and
$\Om^{0+1}(X;\rr) \cap \coker L(A,\rho,t) =\HH^0_A (X;su(3))$
for all $([A],t) \in U^r,$
then $U^{r}$ is a smooth 1-manifold.

If $U \subset W_{\rho} $ is open and regular, then
\begin{enumerate}
\item[(i)]  $W_{\rho}^{*}\cap U $ and $W_\rho^r \cap U$ are both
smooth 1-manifolds without boundary.
\item[(ii)] Each bifurcation point in $U$ is the limit of exactly one
noncompact endpoint of $W_{\rho}^{*},$ i.e.,
$ J = \left( \overline{W_{\rho}^{*}} \setminus W_{\rho}^{*} \right)\cap
U$.
\end{enumerate}

\end{lem}
\begin{pf}
The first statement and (i)
follow from condition (iv) of Definition \ref{p m reg}
using standard Kuranishi arguments.
The proof of (ii) is given  below.

Fix a bifurcation point, which we assume,  for simplicity of notation,
to be of the form  $([A],0)$.
For some neighborhood $U \subset \BB \times [-1,1]$,
$W_{\rho}\cap U$ is the quotient by the gauge group of
the zero set of the map
$$Q:X_{A} \times [-1,1] \lto  \Om^{1}(X;su(3))$$
given by
$ Q(A+a,t) = \Pi_A \zeta_{\rho(t)}(A+a).$

The linearization of $Q$ at $(A,0)$ is an elliptic
Fredholm operator with index one
$$DQ_{(A,0)}: \Om^1(X;su(3)) \oplus \RR \lto  \KK_A $$
and $DQ_{(A,0)}(a,\tau) = \Pi_A L(A,\rho,0)(0,a,\tau).$
Fix a nontrivial $v \in \Om^{1}(X;\rr) \oplus \RR$
in the kernel of $DQ_{(A,0)}.$
Then
$\ker DQ_{(A,0)} = \Span \{v\} \oplus \HH^{1}_{A,\rho_0}(X; \rrp)$
and $\coker DQ_{(A,0)} = \HH^{1}_{A,\rho_0}(X;\rrp)$.

We summarize the Kuranishi model in this situation.  There is a
function $$\phi : \ker DQ_{(A,0)} \lto (\ker DQ_{(A,0)} )^{\perp}$$
and a
neighborhood $U\subset \ker DQ_{(A,0)}$ of zero such that $Q$
restricted to the graph of $\left. \phi \right|_U$ takes values in
$\coker DQ_{(A,0)}$.
Let $\phi_1$ and $\phi_2$ be the $\Om^{1}(X;su(3))$ and $\RR$ components
of $\phi$
and  define the map $\psi : \ker DQ_{(A,0)} \lto \KK_A$
by setting $\psi (a,\tau) = a + \phi_1(a,\tau).$
Now for $s \in \RR$, define
$\Psi_s:\HH^1 _{A,\rho_0}(X;\rr^\perp) \lto X_A$
by setting $\Psi_s (x)= A+\psi (sv+x)$.
Set $\CS_s(A) = \CS(A)+\rho(t_s)(A),$
where $t_s = \phi_2(sv).$
Observe that $t_0 =0.$
Then for all $s$,
$$
Q\circ \Psi_{s} = - 4 \pi^2  \nabla \left( \CS_s \circ \Psi_{s} \right),$$
a family of gradient vector fields of $U(1)$ invariant functions on
$\HH^1 _{A,\rho_0}(X;\rrp)\cong \CC$.

For small $s,$ the path of orbits  $([ \Psi_s(0)],t_s)$
parameterizes $W^{r}_{\rho}$ near
$([A],0)$.  At the origin in $\HH^1 _{A,\rho_0}(X;\rrp)$,
the Hessian of $\CS_s \circ \Psi_s$
is $\la_s \Id$, where $\la_s$ is
the eigenvalue referred to in condition (vi) of Definition \ref{p m reg}.
The proof now reduces to the parameterized Morse Lemma.
See the proof of Theorem 12 in \cite{herald3} for a similar argument.
\end{pf}

Our proof of regularity will involve considering the irreducible and
reducible universal zero
sets
$$ Z^{*}= \left\{ ([A],h)\in \BB^{*}\times \ps \mid \zeta_{h}(A) =0 \right\}$$
and
$$ Z^{r}= \left\{([A],h)\in \BB^{r}\times \ps \mid \zeta_{h}(A) =0 \right\}.$$
Within $Z^{r}$ lies a subset which we hope to avoid when choosing
perturbations, namely, the union over all positive integers $k$ of
$$
Z^{r}_{k} = \left\{ ([A], h)\in Z^{r} \mid  \dim_\CC \ker 
\left( \left. K(A,h) \right|_{\Om^1(X;\rrp)} \right)=k
\right\}.$$

\begin{prop}\label{smooth Z}  The sets
$Z^{*}$ and $Z^r$ are submanifolds of 
$\BB^{*}\times \ps$ and $\BB^{r}\times \ps$, respectively. 
For each $k,$ $Z^r_k$
is a submanifold of $Z^r.$
\end{prop}
\begin{pf}  Fix $([A_{0}],h_{0})\in Z^{*}$.
Consider the map $P:X_{A_0}\times \ps \lto \KK_{A_0}$
given by $P(A,h) = \Pi_{A_0} \zeta_{h}(A)$.
The first partial derivative 
$\frac{\partial P}{\partial a}(A_0,h_{0})$ 
is Fredholm with cokernel $\HH^{1}_{A_0,h_{0}}(X;su(3))$, but, since
$\Ga$  is abundant for $([A_{0}],h_{0})$, the image of
$\frac{\partial P}{\partial h} (A_0,h_{0})$ is a subspace which
orthogonally projects onto this cokernel.  Therefore $P$ is a
submersion at $(A_0,h_{0})$.  The implicit function theorem now proves
that the preimage $P^{-1}(0) \subset X_{A_{0}}\times \ps$ is 
smooth near $(A_{0}, h_{0})$, and hence $Z^{*}$ is smooth near
$([A_{0}],h_{0})$.

To show smoothness of $Z^{r}$, apply the same argument to
the map
$P^{r}:X_{A_0 }^{r} \times \ps \lto \KK_{A_0}\cap \Om^1(X;\rr),$
which is the restriction of the map $P$ to the reducible slice
 $X_{A_{0}} ^{r}=\{A_{0}+a \mid a \in \KK_{A_0} \cap \Om^1(X;\rr)\}$.
That $P^r$ takes values  in $\KK_{A_0}\cap \Om^1(X;\rr) $ follows from 
$\stab A_0$ equivariance. 

Next we treat the third case.  Suppose that $([A_{0}],h_{0})\in
Z^{r}_{k}$.
Define 
$$\la_0 = \min \{ |\la| \neq 0 \mid \la \in \mbox{Spec} (K_{A_0}) \}.$$
Choose a neighborhood $U\times V \subset  X_{A_{0}} ^{r}\times \ps$ of
$([A_{0}],h_{0})$  such that for $(A,h) \in U \times V,$
the operator $K(A,h)$ has no  eigenvalue $\la$ with
$\frac{\la_{0}}{3} < |\la|< \frac{2\la_{0}}{3}$.
Form the small eigenspace bundle, which is the
complex vector bundle $E$ over $U \times V$ with fiber
$E_{(A,h)}$ equal to
$$\Span \left\{ u \in \Om^{0+1}(X;\rr^\perp) \mid K(A,h)(u) = \la u
\hbox{ where } |\la| < \frac{\la_0}{3} \right\}.$$
Let $\sym E$ be the associated fiber  bundle of
symmetric, $\stab (A_0)$ invariant (hence Hermitian)
bilinear forms on $E$,  and for each $k=1,\ldots, \dim_{\CC}
\HH^{1}_{A_{0},h}(X;\rrp)$,
let $\sym_{k} E$ be the
subbundle consisting of those bilinear forms with complex rank less
than or equal to
$\dim_{\CC}
\HH^{1}_{A_{0},h}(X;\rrp)-k$.
Notice that $\sym_{k} E$ has codimension $k^{2}$ in $\sym E$.

Define $\overline{K}(A,h):E_{(A,h)} \lto E_{(A,h)}$
to be the restriction of $K(A,h)$
to $E_{(A,h)}$ composed with the orthogonal projection to $E_{(A,h)},$
and use this to construct the section
$$R:U \times V \lto \sym E \oplus \left(\KK_{A_0}\cap \Om^1(X;\rr) \right)$$
given by
$R(A,h)= \left( \overline{K}(A,h), P^r (A,h) \right).$
Then
$Z_{k}^{r} = R^{-1}\left( \sym_{k} E \oplus 0 \right).$
Now we claim that
$R$ is a submersion at $(A_0,h_0)$.
Since the linearization of $R$ in the first variable has cokernel $T_{0}\sym
E_{(A_0,h_0)} \oplus \HH^{1}_{A_0}(X;\rr),$  it suffices
to show
that the linearization in the other variable,
composed with projection to this cokernel,  is onto.
This is the map
$T_0 \ps \lto \sym \HH^1_{A_0}(X;\rr^\perp) \oplus
\HH^1_{A_0}(X;\rr)$ given by
$$ \de h \mapsto \left(- 4 \pi^2 \hess \de h(A_0),
\; \Pi'_1 \, \nabla \de h(A_0) \right),$$
where
$\Pi'_1$ is the projection onto
${\HH^{1}_{A_0}(X;\rr)}.$
But surjectivity of this map follows since ${\Ga}$ is abundant
for $([A_{0}],h_{0})$.\end{pf}

We are finally ready to prove the regularity theorem for the moduli
space and the parameterized moduli space.   For
$h_{-1},h_{1}\in {\ps}$, let $C^{1}([-1,1], \ps;
h_{-1},h_{1})$ denote the set of $C^{1}$ curves $\rho:[-1,1] \lto
\ps$ with $\rho (\pm 1) = h_{\pm}$.

\begin{thm} \label{regularity}  There exists a Baire
set ${\ps}' \subset \ps$ such that $h\in {\ps}'$
implies $\MM^{*}_{h} \cup \MM^{r}_{h}$ is regular.   For any
$h_{-1}, h_{1}\in {\ps}'$, the set of $\rho \in C^{1}([-1,1],
\ps;
h_{-1},h_{1})$ for which $W_{\rho}$ is regular is Baire.
\end{thm}
\begin{pf}
  The projections
from $Z^{*}$, $Z^{r}$, and $Z^{r}_{k}$ to $\ps$
are Fredholm of index $0, 0, $ and
$-k^{2}$, respectively.  The first two index calculations simply follow
from the
self-adjointness of the partial derivatives in the connection
variable of the maps $P$ and $P^{r}$.  The third follows
easily from the second.
The rest of the argument is a standard 
application of the Sard-Smale
theorem and transversality (see \cite{donaldson kronheimer},
Section 4.3.2).
\end{pf}

\newpage
\section{Orientations and Spectral Flow}
\label{orientations and spectral flow}

In this section, we introduce orientations on the parameterized 
moduli space and relate them to the
spectral flow of the family of
operators $K(A,h)$ from the previous section.  
We use the index bundle of the family $L$
to orient $W^*_\rho$ and $W^r_\rho.$

The basic idea is a familiar one, used not only in 3-dimensional gauge
theory by Taubes (see \cite{taubes}), but also in 4-manifold gauge theory.
In fact, if $W_\rho$ were generically a cobordism,
then Taubes' approach to defining an invariant would work equally
well for $SU(3)$. But $W_\rho$ is {\it not} generically a cobordism,
as explained in
Lemma \ref{regular cobordism},
and a relationship between the orientations on $W^*_\rho$ and $W^r_\rho$
near a bifurcation point is provided by Theorem \ref{end orientations}.

\subsection{Orientations}

Suppose that $\ps$ is fixed as in 
the previous section and
consider the  family of index one Fredholm operators
$$L:\AA\times C^1([-1,1], \ps) \times [-1,1] \lto
\Fred^{1}(\Om^{0+1}(X;su(3)) \oplus \RR, \Om^{0+1}(X;su(3)))$$
introduced in subsection \ref{regularity theorems}.
The dimension of the kernel of $L(A,\rho,t)$ is not continuous in $(A,\rho,t)$,
so $\ker L$ does not form a vector bundle over $\AA \times
C^1([-1,1], \ps) \times [-1,1]$.  
Instead, we consider the index bundle of $L,$ which is
the element in the $K$-theory of $\AA \times C^1([-1,1], \ps) \times [-1,1]$
defined by $\ind L = [\ker L] - [\coker L],$ a virtual bundle of dimension one.

Given vector spaces $E$ and $F$ of dimensions $n$ and $m$, an
orientation on $[E]-[F]$ is an
orientation on the real line $$\det ([E] - [F])= \Lambda^{n}E \otimes
\left(\Lambda^{m}F\right)^{*}.$$
For example, if $\{ e_{1},\ldots,e_{n}\}$ and
$\{f_{1}, \ldots, f_{m}  \}$ are bases for $E$ and $F,$
then the element
$(e_1 \wedge \cdots \wedge e_n) \otimes (f_1 \wedge \cdots \wedge f_m)^*$
specifies an orientation for $[E] - [F].$
More generally, if $E$ and $F$ are vector bundles,
then an orientation on the element $[E] - [F]$ of $K$-theory
is an orientation of the line bundle 
$\Lambda^{n}E \otimes \left( \Lambda^{m}F \right)^{*}.$

Clearly, $\ind L$ is orientable since the parameter space is contractable.  
The virtual fiber at $(\th, 0,0)$ is 
$[\HH^0_\th(X;su(3)) \oplus \RR] -[\HH^0_\th(X;su(3))],$
and our convention for orienting $\ind L$ is to propagate the canonical orientation
at $(\th,0,0)$ given by
\begin{equation}
\label{basis} 
(v_{1} \wedge \cdots \wedge v_{8} \wedge w) \otimes 
(v_{1} \wedge \cdots \wedge v_{8})^*,
\end{equation}
where $\{v_1,\ldots, v_8\}$ is a basis for $su(3) = \HH^0_\th(X;su(3))$
and $w$ is a tangent vector to $[-1,1]$ at $t=0$ pointing 
in the positive direction.

Suppose that $\rho \in C^{1}([-1,1],\ps)$ and $W_{\rho}$
is regular.  Then $W_{\rho}^{*}$ inherits an orientation because
of the natural identification
$T_{([A],t)}W^{*}_{\rho} \cong \ker L(A,\rho, t)$.  
There is also an induced orientation for $W^{r}_{\rho},$ 
but this is less obvious.
First, suppose $([A],t)\in W_{\rho}^{r}$ is not a bifurcation point.
An orientation 
is given by declaring that  a nontrivial vector 
$v \in T_{([A],t)}W^{r}_{\rho}$ is 
positively oriented if
the element $(u \wedge v) \otimes u^* \in \det \ind L (A, \rho, t)$
agrees with the orientation of $\ind L$
for any $u\in u(1) \cong \HH^{0}_{A}(X;su(3)).$

Now suppose that $([A],t) \in W^r_\rho$ is a bifurcation point.
The dimension of $\ker L$ and $\coker L$ both jump by two at $(A, \rho, t)$,
but we obtain an
orientation consistent with the one above by requiring that
$(u \wedge x \wedge y \wedge v) \otimes (u \wedge x \wedge y)^*$
agree with the given orientation on $\ind L(A,\rho,t),$
where
$\{x,y \}$ is a basis for 
$\HH^{1}_{A,h}(X;\rrp)$, the new part of the
kernel (and cokernel) of $L$ at $(A,\rho,t)$.

\subsection{Spectral flow}

In analogy with  Taubes' gauge theoretic description of the Casson 
invariant, our formula will involve counting irreducible perturbed 
flat orbits with sign according to their spectral flow.  
We adopt the following convention for computing
the spectral flow.

\begin{defn}\label{spectral flow} 
Suppose $\UU$ is a real, infinite dimensional, separable Hilbert space
and $K:[0,1] \lto \SAFred(\UU)$ is a 
continuously differentiable family of 
self-adjoint Fredholm operators with discrete spectrum
on $\UU.$ 
Note that the eigenvalues of $K_t$ vary continuously differentiably.
Choose $\de$ such that
$$ 0 < \de   < \inf \{|\la| \neq 0 \mid \la \in \Spec K_0 \cup \Spec K_1 \}.$$ 
The {\bf spectral flow}
along $K_t$ from $K_0$ to $K_1,$  denoted $\SF(K_{0},K_{1}),$ 
is the intersection number, in $[0,1] \times \RR,$
of the graphs of the eigenvalues of $K_{t}$, counted with 
multiplicities, 
with the line segment from $(0,-\de)$ to $(1,\de)$.  
It is a homotopy invariant of the path $K_{t}$ relative to its
endpoints.
\end{defn}
Note that with this convention for counting zero modes,
$$ \SF(K_{0},K_{1})+\SF(K_{1},K_{2})=\SF(K_{0},K_{2}) - \dim \ker K_{1}.$$
We are primarily interested in the spectral flow of the
operator $K(A,h)$ from subsection 
\ref{p flat}.
Completing $\Om^{0+1}(X;su(3))$ in the $L^2$ norm,
we regard
$K(A,h)$ as a family of self-adjoint Fredholm
operators on $\Om^{0+1}(X;su(3))$ with dense domain
the space of $L^{2}_{1}$  forms,
$$K:\AA\times \ps \lto \SAFred\left(\Om^{0+1}(X;su(3))\right).$$

Define
$\deg:\GG \lto \ZZ$ by setting $\deg g = \deg g',$
where $g':X \lto SU(2)$ is a map homotopic to $g$.  
That $\deg g$ is well-defined follows from the next proposition,
which can be proved by noting that $SU(n)$ is homotopy
equivalent to a CW-complex with 3-skeleton $S^3$ and the next
lowest cell in dimension 5.

\begin{prop} \label{degree} 
Fix $n > 2$ and consider the standard inclusion 
$i:SU(2) \subset SU(n).$ 
\begin{enumerate}
\item[(i)] If
$g \in C^\infty(X,SU(n)),$ 
then there exists
$g':X \lto SU(2)$ with $i \circ g' \simeq g.$
\item[(ii)] If $g_0,g_1 \in C^\infty(X,SU(2))$ with $i \circ g_0 \simeq i \circ g_1,$
then $g_0 \simeq g_1.$
\end{enumerate}
\end{prop}

Proposition \ref{degree} gives the following formula for the
spectral flow between two gauge equivalent connections.

\begin{prop}  
\begin{enumerate}
\item[(i)] The spectral flow of $K(A,h)$ along a path $(A_t,h_t)$
is independent of the path connecting  $(A_{0},h_{0})$ to
$(A_{1},h_{1}).$ 
\item[(ii)]  The spectral flow of $K$ from $(A,h)$ to $(gA,h)$ equals $12 \deg g-\dim \ker K(A,h)$.
\end{enumerate}
\label{spec flow and deg}
\end{prop}
\begin{pf}  Part (i) follows since $\AA\times \ps$ is contractable.
Part (ii) follows by an index computation, the point being that
spectral flow around a closed path in $\AA$ equals the index of the
self-duality operator on $SU(3)$ connections over $X \times S^1.$ 
Details can be found in \cite{kkr}.
\end{pf}

\noindent{\it Remark.} Suppose $A\in \AA^{r}.$ By applying
a gauge transformation, we can assume that $A \in \AA_{S(U(2) \times U(1))}.$
Consider now the standard decomposition
of the Lie algebra
$su(3)= \rr \oplus \rrp$ given by
the action of $\stab A = U(1)$ 
and split the operator $K(A,h)$
accordingly. Because $\th$ and $A$ can be
connected by a path 
in $\AA_{S(U(2)\times U(1))},$
the spectral flow of $K$ from $(\th, 0)$ to 
$(A,h)$ splits as
$$\SF(\th,A) = \SF_{\rr}(\th,A) 
+ \SF_{\rrp}(\th,A).$$
Notice that $U(1)$ equivariance of $K(A,h)$ implies that
$\SF_{\rrp}(\th, A)$ is
divisible by two. Using part (ii) of the previous proposition
and the well-known, analogous result (for $su(2)$)
that $\SF_\rr(A,g A) = 8 \deg g - \dim \ker K_A \vert_{\Om^{0+1}(X;\rr)},$ 
we see that 
$$\SF_\rrp(A,g A) = 4 \deg g - \dim \ker K_A \vert_{\Om^{0+1}(X;\rrp)}.$$

\subsection{The relationship between orientations and spectral flow}

There is a fundamental relationship between the orientation of the
one dimensional virtual bundle $\ind L$ and the spectral flow of 
$K(A,h)$.   We describe it next, in some generality.

Suppose that $\UU$ is an infinite dimensional, separable  Hilbert space
and that $Z$ is a connected, simply connected parameter space.
Let $$K:Z \lto \SAFred(\UU) $$ be a parameterized family of
self-adjoint Fredholm operators on $\UU$    
and $v:Z \lto \UU$ be a continuous map.  

Define $L_z:\UU \oplus \RR \lto \UU$ by $L_z(u,\tau) =
K_z(u)+ \tau v_z$ for $(u,\tau) \in \UU \oplus \RR.$
Clearly $L_z \in \Fred^1(\UU \oplus \RR,\UU)$.
For any $z \in Z,$ let $\Pi_z:\ker L_z \lto \RR$ be the
projection onto $\RR$
and $\Pi_{\ker L_z}: \UU \oplus \RR \lto \ker L_z$ be the
projection onto the $\ker L_z.$ 

Suppose that  $z_{0}\in Z$ is a fixed base point and $v_{z_{0}} =0$.
Choose an orientation $\OO$ for $\ind L$ by the convention in
equation (\ref{basis}), and let $\OO_{z}$ denote the induced orientation on
$\ind L_{z}$.
If $L_z$ is surjective, then $\OO_z$ gives an
orientation of $\ker L_z.$ Notice that
whenever $K_z$ is an isomorphism, $\ker L_z$ is spanned by
$( -K_z^{-1}(v_z), 1).$
In this case, and more generally when $v_{z}\perp \ker K_{z}$, 
the  spectral flow of $K_z$ allows us to compare $\OO_z$
with another natural
orientation on $\ker L_z$.

\begin{prop}\label{abstract version}
Suppose $v_{z_{1}} \perp \ker K_{z_{1}}$.  Then if $\{u_{1}, \ldots, 
u_{k}\}$ is a basis for $\ker K_{z_{1}}$, $\{ u_{1}, \ldots, u_{k}, 
\left(-K_{z_{1}}^{-1}(v_{z_{1}}), 1\right) \}$ is a basis for $\ker 
L_{z_{1}}$.  Furthermore, the orientation on $\ind L_{z_{1}}$ 
agrees with 
$$ (-1)^{\SF (K_{z_{0}}, K_{z_{1}})}\left( u_{1}\wedge \ldots \wedge 
u_{k} \wedge \left(-K_{z_{1}}^{-1}(v_{z_{1}}),1 \right)\right)  \otimes 
(u_{1}^{*}\wedge \ldots \wedge u_{k})^{*}.$$
\end{prop}
\noindent{\em Remark.}  $K_{z_{1}}$ is an isomorphism
if and only if $\Pi_z$ is an isomorphism, and then the 
proposition states that the orientation on $\ker L_{z_{1}}$ is 
$(-1)^{\SF (K_{z_{0}}, K_{z_{1}})} \Pi_{z_{1}}^{*} \OO_{\RR}$.  
\begin{pf}  The first claim is obvious. 
The proof of the second goes as follows.  
Connect $z_{0}$ to $z_{1}$ by a path $z_{t}$.  By \cite{koschorke},
$K_{z_{t}}$ is homotopic relative its endpoints to a
path $K_{t}$ in $\SAFred(\UU)$ so that there is a finite set
$\{t_1,\ldots, t_k\} \subset (0,1)$ such that
$$\dim \ker K_t = \left\{
\begin{array}{ll} 1 & \hbox{if } t \in \{t_1,\ldots, t_k\} \\
					  0 & \hbox{otherwise.}
\end{array} \right.$$
We can further assume that any 
eigenvalue of $K_{t} $ which crosses zero does so transversely.

Similarly, $v_{z_{t}}$ can be homotoped relative to its endpoints to
a path $v_{t}$ in $\UU$ such that the path $L_{t}$ in
$\Fred^{1}(\UU\oplus \RR,\UU)$ defined by $L_{t}(u,\tau) = K_{t}(u) +
\tau v_{t}$ is surjective for all $t\in (0,1]$.
Let $\OO_{t}$ be the orientation on $\ker L_{t}$ coming from
$\OO_{0}=\OO_{z_{0}}$.

Fix a $t_{j}$ with $\ker K_{t_{j}}$ nontrivial.
For $t\in (t_{j }-\de, t_{j} + \de),$ 
let $\la_{t}$
be the eigenvalue of
$K_{t}$ which crosses zero when $t=t_{j}$.  
Choose $u_t$ to be a unit eigenvector with eigenvalue $\la_t$ 
so that $K_{t}(u_{t})= \la_{t} \cdot u_{t}$.

For $t\in (t_{j}-\de, t_{j}+ \de)$, we have an orthogonal
decomposition of $\UU$ into
$\UU'_t \oplus \UU''_t$ where $\UU''_t = \Span \{ u_{t} \}$ and
$\UU'_t$ is its orthogonal complement.
Set $a_t = \langle u_t, v_t \rangle, \;
v'_t = v_t - a_t u_t$ and 
$K'_t(w) = K(w) - \la_t \langle u_t,w \rangle u_t$ for $w \in \UU.$
Note that $K'_t$ is invertible on $\UU'_t$ 
and set $w_{t}=a_{t}u_{t} +
\la_{t}(K'_{t})^{-1}v'_{t}$.  
The vector $(w_{t},-\la_t)$ spans $\ker L_{t}$ for $t\in (t_{j}-\ep, t_{j}+ \ep)$.
Since the inner
product $\langle (w_{t}, -\la_t), (0,1) \rangle$ changes sign at $t_{j}$, it
follows that the orientation of $\Pi_t^* \OO_{\RR}$ changes relative to $\OO_t$
at $t_j.$ Such a change occurs for each $t_j,$ which is where
$\SF(K_0, K_t)$ changes by $\pm 1.$

This proves the second claim in case $\ker K_{z_{1}}$ is 
trivial.  For the general case, we may assume that
all the eigenvalues 
of $K_{t}$ which approach zero as $t \to 1^-$ 
are negative for $t$ near 1. 
This implies that $\SF(K_{0}, K_{1})=\SF(K_{0}, K_{t})$ 
for $t\in (1-\de, 1)$.  We then claim 
that the orientation given by  
$\Pi_{\ker L_{t}} \left(-K_{z_{t}}^{-1}(v_{z_{t}}),1 \right)$
propagates to 
$\left( u_{1}\wedge \ldots \wedge 
u_{k} \wedge \left(-K_{z_{1}}^{-1}(v_{z_{1}}),1 \right)\right) \otimes 
(u_{1}\wedge \ldots \wedge u_{k})^{*}$. 

Recall our convention for propagating the orientation of
$\ind L$ across a point where the $\dim \ker L$ jumps.  
The orientation given by  
$\Pi_{\ker L_{t}} \left(-K_{z_{1}}^{-1}(v_{z_{1}}),1 \right)$
propagates to 
$$\left( \left(-K_{z_{1}}^{-1}(v_{z_{1}}),1 \right)\wedge 
u_{1}\wedge \ldots \wedge u_{k}\right) \otimes 
\left( \widetilde{L}_{t}(u_{1})\wedge \ldots \wedge \widetilde{L}_{t}(u_{k}) \right)^{*},$$ 
where $\widetilde{L}_{t} = \Pi_{\coker L_1} \circ L_t$. 
Since $\widetilde{L}_{t}$ is 
negative definite on $\Span \{u_{1}, \ldots, u_k\}$, it follows that 
$\widetilde{L}_{t} (u_{1})\wedge \ldots \wedge \widetilde{L}_{t} (u_{k})$
is proportional to $(-1)^{k} u_{1}\wedge \ldots \wedge u_{k}$.  
Permuting the $\left(-K_{z_{1}}^{-1}(v_{z_{1}}),1 \right)$ factor 
past all the $u_{i}$'s  
introduces another $(-1)^{k}$ which cancels with the first.
\end{pf}

Applying Proposition \ref{abstract version} to the
oriented strata in a regular moduli space gives the following
corollary. 
\begin{cor} \label{orientations sf}
Assume that $\rho:[-1,1] \lto \ps$ is a path of perturbations 
such that $\MM_{\rho(+1)},$ $ \MM_{\rho(-1)}$, and $W_{\rho}$ are all regular.  
Then $-1$ and $+1$ are regular values of the projections from 
$W_{\rho}^{*}$ and $W_{\rho}^{r}$ to $[-1,1].$
Suppose $\varepsilon = \pm 1$ and
$([A], \varepsilon) \in 
\MM^{*}_{\rho(\varepsilon)} \cup \MM^{r} _{\rho(\varepsilon)}$,
and set $s = \SF(K_{\th, 0}, K_{A,\rho(\varepsilon)})$.
Then the boundary 
orientation of $W_{\rho}^{*}$ or $W_{\rho}^{r}$ at $([A],\varepsilon)$ 
equals 
$(-1)^s$ if $\varepsilon = 1$ and it equals $-(-1)^s$
if $ \varepsilon = -1.$
\end{cor}
\begin{pf} 
Note that the boundary orientation at $([A],\varepsilon)$ is positive  if 
and only if the orientation on the 1-dimensional stratum of $W_{\rho}$ 
at $([A],\varepsilon)$ agrees with $\varepsilon \Pi^{*}\OO_{\RR}$.  
In the irreducible case, $K(A,\rho(\varepsilon))$ is an isomorphism, so the 
remark following Proposition \ref{abstract version} proves the claim.

The reducible case 
also follows by a direct application of Proposition \ref{abstract version},
letting $-4 \pi^2 \left. \frac{\partial}{\partial t} 
\nabla \rho_{t}(A)\right|_{t=\varepsilon}$ play the role of the 
$v_{z_{1}}$ for the operator $L(A,\rho,t)$ and observing that this vector
is orthogonal to $\ker K(A,\rho(\varepsilon))= 
\HH^{0}_{A}(X;su(3))$.  
\end{pf}
  
\subsection{Orientations near a bifurcation point}
In this subsection, we identify the boundary orientation of
a bifurcation point with the oriented $\rrp$ spectral flow of $K(A,h)$
along $W_{\rho}^{r}$ across this point.
The precise relationship is
given in Lemma \ref{calc}. This is the crucial
observation needed for Theorem \ref{end orientations},
which is used in section \ref{the invariant}
to show that our invariant 
is well-defined.

Consider  the operator 
$ L(A,\rho,t):\Om^{0+1}(X;su(3)) \oplus \RR \lto \Om^{0+1}(X;su(3))$
for a fixed $\rho \in C^1([-1,1], \ps)$
such that $W_{\rho}$, $\MM_{\rho(-1)}$, and $\MM_{\rho(+1)}$ are regular.
Suppose that $W_\rho$ has a bifurcation point, which we 
take to be  
$([A],0)$ for simplicity of notation.  Assume that $A\in 
\AA_{S(U(2)\times U(1))}$ is a 
 representative of the orbit $[A]$. 
Choose   a  covariantly  constant, diagonal $su(3)$-valued
0-form 
$${u}=\left( \begin{array}{ccc} i /3&&0\\
                           &i /3\\
                           0&&-2i/ 3 \end{array} \right)\in 
                           \HH^{0}_{A}(X;su(3)).$$
Then the complex structure $J$ on
$\Om^{0+1}(X;\rrp)$ is given by $\exp
 ({\pi  {u}/2})\in \stab A $  acting by conjugation, 
i.e.,
$J x = [u,x]$ for $x \in \Om^{0+1}(X;\rrp).$
Choose a nonzero  $x \in \HH^1_{A,\rho_0}(X;\rrp) $ and set $y=Jx$.

Let
$v\in \Om^{1}(X;\rr)\oplus \RR,$ be an element of $\ker L(A,\rho, 0)$ 
such that
\begin{equation} \label{orient1} 
(u \wedge x \wedge y \wedge v) \otimes (u \wedge x \wedge y)^*
\end{equation}
is the orientation for $\ind L$ at $(A,\rho,0).$  In other words, $v$ is 
an oriented tangent vector for $W^{r}_{\rho}$.

Solutions to the equation $\zeta_{\rho(t)}(A')=0$ near $(A,0)$ in 
$X_{A}\times [-1,1]$ take the form $(A',t)=(A+sx + o(s^{2}), 
o(s^{2}))$, $s>0$, up to the action of $\stab A$.  For such a nearby 
solution, $x$ projects nontrivially into the 1-dimensional 
kernel of $L(A', \rho, t)$ (this 
follows from Lemma \ref{regular cobordism}) and its image, 
thought of as a tangent vector to $W^{*}_{\rho}$, points away from 
the endpoint.

We shall now compare the orientation of $\ind L$ at $(A',\rho,t)$ 
with that given by $x.$ To do so, we consider
$$L'(\xi,a,\tau) =  \left. 
{\tfrac{\partial}{\partial s}}  
L(A + s x,\rho, 0)(\xi,a,\tau) \right|_{s=0},$$
where the map on the right is restricted to $\ker L(A,\rho,0)$ 
and then projected onto $\coker L(A,\rho,0)$.

One can check that  $\dim \ker L' = 1,$   and so 
the orientation on $\ker L(A',\rho, t)$
points in the direction of $\ep x,$  where $\ep = \pm 1$ is such that
\begin{equation}\label{orient2} 
(\ep x \wedge u \wedge v \wedge y) \otimes (L'(u) \wedge L'(v) \wedge L'(y))^*
\end{equation}
is the orientation  for $\ind L$ at $(A,\rho,0)$. 
The following lemma is the key step in proving 
Theorem \ref{end orientations} because it
identifies this $\ep$ in terms of 
the $\rrp$ spectral flow of $K(A,h)$.
 
\begin{lem} \label{calc} Suppose that $x,y,u$ are chosen as above. 
Denote by $L_1'$ the composition of $L'$ with
the projection onto $\Om^1(X;su(3)).$ Then
\begin{enumerate}
\item[(i)] $L' ({u}) = -y$ and $L' (y) = -{u}$, and
\item [(ii)] $L_1'(v) = D_{v} ( *d_{A',\rho_{t}} (x))|_{(A',t)=(A,0)}.$ 
\end{enumerate}
\end{lem}
\noindent{\em Remark.}  Recall from section 2
that $*d_{A,h} = *d_A - 4 \pi^2 \hess h(A).$  The notation $D_{v}$ 
means the derivative as $(A,t) $ is varied with tangent vector $v$.
\begin{pf}
First we compute that
\begin{eqnarray*}
L'(u) &=& \left.  
{\tfrac{\partial}{ \partial s}}   L(A+sx,0)(u,0,0) \right|_{s=0}\\
&=& \left.  {\tfrac{\partial}{ \partial s}} d_{A+sx} (u) \right|_{s=0} \\
&=& [x,u] = - [u,x] = -Jx = -y.
\end{eqnarray*}
A similar computation yields $L'(y)=-u,$ 
and these together prove (i).

Claim (ii) follows by commuting mixed partials as follows.  Let 
$(a,\tau)$ denote the components of $v$ in $\Om^{1}(X;\rr) \oplus 
\RR$.  
\begin{eqnarray*}
L'_1(a,\tau) &=& \left. {\tfrac{\partial}{\partial s}} 
\left( * d_{A+sx} (a) - 4 \pi^2 \hess \rho_0(A+sx)(a) - 4 \pi^2 \tau \left. 
{\tfrac{\partial}{\partial t}} 
\nabla \rho_t (A+sx) \right|_{t=0} \right) \right|_{s=0} \\
&=& \left. {\tfrac{\partial}{ \partial s} \frac{\partial}{ \partial r}}
\left( *F(A+sx+ra)- 4 \pi^2 \nabla \rho_{r \tau} (A+sx+ra) \right) 
\right|_{(r,s)=(0,0)} \\
&=& \left. D_{(a,\tau)}\left( \left. 
{\tfrac{\partial}{\partial s}}
\left( *F(A'+sx) - 4 \pi^2 \nabla \rho_t(A'+sx) \right) \right|_{s=0} \right) \right|_{(A',t)=(A,0)}\\
&=& \left. D_{(a,\tau)} (* d_{A',\rho_t}(x))\right|_{(A',t)=(A,0)}
\end{eqnarray*}
This completes the proof of (ii). \end{pf}

Using part (i) of Lemma \ref{calc} and comparing the
two orientations for $\ind L$
at $([A],0)$ given in equations (\ref{orient1}) and $(\ref{orient2}),$
we see that $\ep$ has the opposite sign of the inner product
$\langle L' (v), x\rangle,$ where $v$ is
the oriented vector tangent to $W^r_\rho$ at $([A],0)$.
From  part (ii) of the lemma, 
it follows that $\langle L' (v), x\rangle $ has
the same sign as the derivative of the path of (multiplicity two)
eigenvalues of
$K(A+ra,\rho_{r\tau})$ which crosses zero at $r=0$.
Figure 2 illustrates what this means in terms of 
spectral flow.  Here $W_\rho^r$ is the dotted line and
$W_\rho^*$ the solid one.
In the diagram on the left, $\SF_\rrp(A_-, A_+) = -2$
and in the one on the right, $\SF_\rrp(A_-, A_+) = 2.$

\begin{figure}[hbt]
\begin{center}
\leavevmode\hbox{%
\epsfxsize=4in
\epsffile{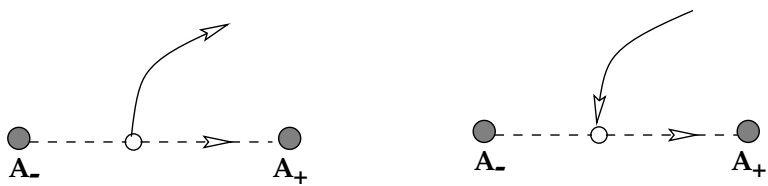}}
\end{center}
\caption{A neighborhood of a bifurcation point in $W_\rho$.}
\label{figure2}
\end{figure}

\noindent
{\it Caution.} The operator $K(A,h)$ on $\Om^{0+1}(X;\rrp)$ is 
equivariant with respect to the action of $\stab A \cong U(1),$ 
thus it is Hermitian with respect to $J$. 
Viewed this way,
the eigenvalue here would have (complex) multiplicity one,
but in order 
to avoid confusion, we regard $K(A,h)$ on (0+1)-forms with values in
either component of the
splitting $su(3) = \rr \oplus \rrp$ as
a {\em real} operator.

Summing over all the bifurcation points results
in the following theorem.

\begin{thm}\label{end orientations}
Let $\rho \in C^{1}([-1,1],\ps)$ 
be a curve with $W_{\rho}, \MM_{\rho(\pm1)}$ regular and suppose
$C$ is a connected component of  $W_{\rho}^{r}.$ 
Define $b(C)$ to be the number of bifurcation points on $C$ counted with
orientation as boundary  points of $\overline{W_{\rho}^{*}}$.
\begin{enumerate}
\item[(i)] If  $\partial C = \emptyset,$ then $b(C)=0.$
\item[(ii)] If $\partial C = ([A_+],\varepsilon_+) \cup -([A_-],\varepsilon_-),$ 
where $\varepsilon_\pm \in \{-1,+1\} $ are not necessarily distinct, 
then  $b(C)=\frac{1}{2}\SF_{\rrp}(A_{-},A_+),$
provided the representatives $A_\pm$ for $[A_\pm]$ are chosen to lie
on the same component of a lift of $C$ to $\AA^r.$
\end{enumerate}
\end{thm}

\begin{pf}
To prove (i), suppose $C$ is a component 
of $W_{\rho}^{r}$ with $\partial C=\emptyset.$ 
Choose a path of connections $A_s$ and perturbations $h_s$
for $s \in [0,1]$ such that $([A_s],h_s)$ parameterizes $C$.
Then $A_1 = g A_0$ for some $g \in \GG.$
Proposition \ref{existence of epsilon} implies that
the entire path $A_s$ of perturbed flat connections
lies within $\ep_0$ of some component $K$ of the space
of flat connections upstairs in $\AA.$
By our choice of $\ep_0,$
that proposition also shows that $A_s$ does not come within
$\ep_0$ of any other component of the flat connections.
But if $A_0$ is within $\ep_0$ of $K,$ then $gA_0 = A_1$
is within $\ep_0$ of $g K$ since we are using the standard
gauge invariant $L^2_1$ metric. Thus $gK=K.$

Now $\CS: \AA \lto \RR$ is constant on connected components of
the space of flat connections, and 
this implies that $g \in \GG_0,$ the connected component of
the identity in $\GG,$ since otherwise $\deg g \neq 0.$
Therefore, using the relationship between spectral flow
and degree described in the remark following
Proposition \ref{spec flow and deg}, we get
$$b(C) = {\tfrac{1}{2}} \SF_{\rrp}(K_{A_{0},h_{0}},
K_{A_{1},h_{1}})
=2( \CS(A_{1}) -\CS(A_{0})) = 0.$$
This proves (i), and part (ii) of the theorem is clear.
\end{pf}

\noindent {\it Example.} We indicate briefly the consequence of
the above theorem
for the situation illustrated in
Figure \ref{figure1}. 
First of all, part (i) of Theorem \ref{end orientations} implies
that the 
$\rrp$ spectral flow around the one closed component equals 0.

Along the other components, which are the two other dotted curves,
the $\rrp$ spectral flow in the oriented direction equals 2
for the component on top and
it equals 4 for the one on bottom.
In other words, the $\rrp$ spectral flow along the bottom component from
{\it left} to {\it right} equals $-4.$

\newpage
\section{The Invariant} \label{the invariant}
In this section, we define the invariant $\la_{SU(3)}(X)$ 
for $X$ an orientable, integral homology 3-sphere. 
Choose an orientation and Riemannian metric on $X,$
as well as a collection $\Ga = \{ \ga_1, \ldots, \ga_n \}$
of embedded solid tori in $X$ satisfying the conlusion of
Proposition \ref{openness}.
Let $\ps$ be the $\ep_0$ neighborhood of 0 in ${\cal F}_\Ga,$
where $\ep_0$ is given by Proposition \ref{existence of epsilon}.
Then choose a perturbation $h \in \ps$ so that $\MM_h$ is regular.
By Proposition \ref{regularity1},
$\MM_h^*$ is a compact 0-manifold. 
We would like to define an invariant of $X$ by counting
the points $[A] \in \MM^*_h$ with sign according to the parity
of the spectral flow of $K.$
This integer, however, depends on the choice of perturbation $h \in \ps$
and in order to obtain a well-defined invariant, 
we must include a correction term
determined from $\MM_h^r.$
 
When the perturbation $h$ is clear from the context, we 
let $\SF(A_0,A_1)$ be an abbreviation for
$\SF(K_{A_0, h}, K_{A_1,h})$.
For $A_0,A_1 \in \AA^r$, the spectral flow splits as
$$\SF(A_0,A_1) = \SF_\rr(A_0,A_1) + \SF_\rrp(A_0,A_1)$$
according to the decomposition $su(3) = \rr \oplus \rrp$.

\begin{thm} \label{indie}
Suppose that $h\in \ps$ with  $\MM_{h}$ is regular.
Pick gauge representatives $A$ for each orbit $[A] \in \MM_h,$
and for each representative of a reducible orbit $[A] \in \MM^r_h,$
choose also a
flat connection $\widehat{A}$ with
$\| A-\widehat A \|_{L^2_1} < \ep_0.$
The quantity
$$ \sum_{[A] \in \MM^*_h} (-1)^{\SF(\th, A)}
- \frac{1}{2}\sum_{[A] \in \MM^{r}_h}(-1)^{\SF(\th, A)}
(\SF_{\rrp}(\th, A) - 4 \CS(\widehat{A})).$$
is
independent of choice of representatives $A$ for $[A]$ in both sums and independent
of the choice of $h.$
\end{thm}
 
\begin{pf}
Note that the existence of $\widehat{A}$ is
guaranteed by Proposition \ref{existence of epsilon}.
We first argue that the quantity is independent of the
representatives $A$ chosen for the orbits $[A] \in \MM_h^r.$
Write
$\la'(h) = \sum_{[A] \in \MM^*_h} (-1)^{\SF(\th, A)}$
for the first sum and $\la''(h) = \frac{1}{2}\sum_{[A] \in \MM^{r}_h}(-1)^{\SF(\th, A)}
(\SF_{\CC^2}(\th, A) - 4 \CS(\widehat{A}))$
for the second.
Part (ii) of Proposition \ref{spec flow and deg} shows that
$\la'(h)$ is  independent of choice of the representatives $A$ for $[A] \in \MM^*_h.$
Also, for a fixed representative $A$ for $[A] \in \MM^r_h,$
if $\widehat{A}$ and $\widehat{A}'$ are both flat connections in an
$\ep_0$ neighborhood of $A,$ then  part (iv) of  Proposition \ref{existence of epsilon}
implies that $\widehat{A}$ and $\widehat{A}'$ lie on the same component of the flat connections,
hence $\CS(\widehat{A}) = \CS(\widehat{A}').$
Now suppose $g \in \GG.$ 
Then, by the remark following Proposition \ref{spec flow and deg},
$\SF_{\rrp}(A,g A) = 4 \deg g$.  
Since  $\CS(g \widehat{A}) = \CS(\widehat{A})+\deg g,$
$\la''(h)$ is also independent of the choice of representatives $A$ for $[A] \in \MM_h^r.$

We now argue that the above quantity is independent of choice of $h.$
Suppose that $h_{-}$ and $h_{+}$ are admissible functions in 
$\ps$  and that $\MM_{h_{\pm}}$ are both regular.  Set 
$\MM_{\pm}=\MM_{h_{\pm}}$ and  connect $h_{-}$ and $h_{+}$ by
a path $\rho: [-1,1] \to \ps $ with $\rho(\pm 1) = h_\pm$ so that 
the parameterized moduli space $W_{\rho}$ is regular.

Compactify the irreducible stratum ${W}_{\rho}^{*}$
by adding bifurcation points and denote the compact, 
oriented 1-manifold with boundary
so obtained by $\overline{W_\rho^*}$.
Of course, the total number of boundary points, counted with 
boundary orientation, equals zero.  Every boundary point which is
not a bifurcation point can be 
identified with a point in the disjoint union
$\MM^*_{-} \cup \MM^*_{+}.$  The orientations of these points
are described by Corollary \ref{orientations sf}, as follows.
For $[A]\in \MM^{*}_{+}$, the boundary orientation 
of $W^*_\rho$ at $([A],+1)$ is $(-1)^{\SF(\th, A)}$, 
while for 
$[A]\in \MM^{*}_{-}$, the boundary orientation 
of $([A],-1)$ at  $W^*_\rho$ is $-(-1)^{\SF(\th,A)}$. 
Therefore  
$\la'({h_+}) - \la'({h_-})$
equals minus the number of bifurcation points counted with orientation as
boundary points of
$\overline{W_\rho^*}.$  

It remains to show that this algebraic sum of bifurcation points 
equals $ \la''({h_+}) - \la''({h_-}).$ 
To prove this, we invoke  Theorem \ref{end orientations}.
By part (i), the closed components of 
$W^{r}_{\rho}$ do not contribute to this sum,
so suppose that $C$ is a component of $W^r_\rho$ and
$\partial C = ([A_{+}],\varepsilon_{+}) \cup -([A_{-}],\varepsilon_{-}),$ where 
$\varepsilon_\pm \in \{-1,1\}.$
Let $b(C)$ be the algebraic sum of bifurcation points on $C.$
Since $\la''(h_+)$ and $\la''(h_-)$ are both independent of the choice
of representatives $A$ for $[A],$
we can choose $A_+$ and $A_-$ to lie on the same component of the lift of $C$ to $\AA^r.$
Thus $\CS(\widehat{A_{+}})= \CS(\widehat{A_{-}})$.
By part (ii) of Theorem \ref{end orientations}, 
$$b(C)={\tfrac 1 2} \SF_{\rrp}(A_{-}, A_+) 
= \tfrac 1 2 \left(\SF_{\rrp}(\th, A_{+}) - \SF_{\rrp}(\th, A_{-}) \right).$$
On the other hand, the contribution to $\la''(h_{+})-\la''(h_{-})$ 
from the endpoints of $C$ is 
$${\tfrac 1 2 } \left[ \varepsilon_{+}(-1)^{\SF(\th, A_{+})}\SF_{\rrp}(\th, 
A_{+}) + \varepsilon_{-}(-1)^{\SF (\th , A_{-})} \SF_{\rrp}(\th, A_{-}) 
\right].$$
It is important to keep in mind that
$\varepsilon_{\pm}$ need not be distinct; several 
possibilities are pictured in Figure 1.
Now the reducible case of Corollary 
\ref{orientations sf} implies that $\varepsilon_{+}=(-1)^{\SF(\th, A_{+})}$
and $\varepsilon_{-} = (-1)^{\SF(\th, A_{-})}$, 
and this completes the proof. \end{pf}

The quantity in Proposition \ref{indie}
is seen to be independent of the choice of metric on $X$ 
by the same argument as was used for Proposition 2.3 of \cite{taubes}.
That it is also independent of the choice of $\Ga$ is an exercise which
we leave for the reader.
\begin{defn}  \label{5.3} 
Suppose that $h\in \ps$
 and that $\MM_{h}$ is regular.  
Define the {\bf SU(3) Casson invariant} by
$$\la_{SU(3)}(X) = \sum_{[A] \in \MM^*_h} (-1)^{\SF(\th, A)}
- \frac{1}{2}\sum_{[A] \in \MM^{r}_h}(-1)^{\SF(\th, A)}
(\SF_{\rrp}(\th, A) - 4 \CS(\widehat{A}) + 2).$$
By Theorem \ref{indie}, this gives a well-defined invariant of integral homology 3-spheres.
\end{defn}

Notice that the last term in the second sum above simply adds a multiple of the
$SU(2)$ Casson invariant. This part of $\la_{SU(3)}(X)$
is independent of $h$ by the argument given in \cite{taubes}.
Therefore, the previous theorem implies that $\la_{SU(3)}(X)$
is independent of $h \in \ps$. 
The following proposition explains why we have chosen
to normalize $\la_{SU(3)}(X)$ this way.
\begin{prop}
\begin{enumerate}
\item[(i)] If $\pi_1 X = 0,$ then $\la_{SU(3)}(X)=0.$
\item[(ii)] $\la_{SU(3)}(-X) = \la_{SU(3)}(X).$
\end{enumerate} 
\end{prop}
 
\begin{pf} Part (i) is obvious. To prove (ii), observe that  
$$\SF_{-X}(A_0,A_1) = - \SF_X(A_0,A_1) -(\dim \ker K_{A_0} + \dim \ker K_{A_1}),$$
where the subscript indicates a choice of orientation on $X.$
This is equally valid for 
$\rrp$ coefficients in case $A_0$ and $A_1$ are reducible.
Applying this to all three spectral flows
appearing in the definition of $\la_{SU(3)}(-X)$ and noting
further that
$\CS_{-X}(\widehat{A}) = - \CS_X(\widehat{A})$  
complete the proof of part (ii).
\end{pf}

\newpage
\section{Existence of Perturbation Curves}
This section is devoted to finding loops in $X$ with certain
properties  
required for our transversality arguments in Section 3.
The basic question is whether the trace of holonomy
can detect
a tangent vector to the flat moduli space.
In terms of
a one-parameter family
$A_t$ of irreducible flat $SU(3)$ connections,
we ask: does there exist an element $\ga \in \pi_1(X)$ such that
$$\left. {\tfrac{d}{dt}} \tr\hol_\ga (A_t) \right|_{t=0} \neq 0?$$
The answer is no if $A_t = g_t\, A_0,$ so 
we must also assume that $A_t$ is not tangent to the  
gauge orbit $\GG A_0.$
In fact, we need this for any path $A_t$ of connections
such that $A_0$ is flat and $A_t$ is flat to first order 
(i.e., $\left. \frac{d}{dt} F_{A_t} \right|_{t=0} =0$). 
An affirmative answer to this question for $SU(2)$ and $SU(3)$ is given in the
first two subsections. The last subsection treats the reducible case, where
second order arguments are required.

\subsection{First order arguments} \label{s6.1}

To start, we introduce some notation.
Given a flat connection $A$ and a based loop $\ell : [0,1] \lto X,$ let
$H_\ell(A) \in SU(3)$ be the holonomy of $A$ around $\ell.$
For $a \in \Om^1(X, su(3)),$ let
$I_\ell(a,A) \in su(3)$ be the integral
$$I_\ell(a,A) =\int_0^1 P_\ell (0,t)^{-1} a_{\ell(t)} P_\ell(0,t) dt,$$
where $P_\ell(0,t)$ is the parallel translation from $0$ to $t$ 
along $\ell$ using the connection $A.$
When $A$ and $a$ are clear from context, we write simply
$H_\ell$ and $I_\ell.$

If $A_t = A +ta + O(t^2),$ by Corollary \ref{derivative of tr hol}, we see that 
$\left. \frac{d}{dt} \tr H_\ell(A_t) \right|_{t=0} 
=\tr(H_\ell(A) I_\ell(a,A)).$
\begin{prop}
Suppose $A$ is a flat $SU(3)$ connection.
If $a \in \HH^1_A(X;su(3))$
is non-zero,  
then there is a
loop $\ell$ so that $I_\ell $ projects non-trivially onto $z(H_\ell),$
the Lie algebra of the centralizer $Z(H_\ell).$
\label{first loop}
\end{prop}
\begin{pf}
Consider 
the differential equation $d_A \xi = a$.
We can solve this equation locally on any 3-ball $B \subset X$ since 
$H^1 _A (B; su(3))=0.$ 
Because $a$ is not exact, there is no global solution.
Thus there exists some loop
$\ell :[0,1] \lto X$ (which can be taken to be embedded)
for which the local solutions do
not match up at the ends.
Hence $a|_{\ell}$ is not exact, 
meaning that its decomposition $a|_{\ell}=a_{h}+ d_{A}b$ into 
harmonic and exact parts has $a_h \neq 0.$
From this point on, we restrict our attention to the
pullback connection $\ell^*(A)$ on the $SU(3)$ bundle
over the circle $S^1 = [0,1]/0 \sim 1$ pulled back via the loop $\ell.$

Note that $a_h$ is Hodge dual
with respect to the metric on the loop 
to a covariantly constant 0-form, and so integrates to
something nonzero in $\HH ^{0}_{\ell^*(A)}(S^1;su(3))$, the Lie algebra of 
$\stab (\ell^*(A))$.
By the fundamental theorem of calculus,
the exact part integrates to 
$H_{\ell}^{-1}b_0 H_{\ell} - b_0$, where $b_0$ denotes the value of
$b$ at the basepoint.  This latter $su(3)$ element is orthogonal to
$\HH^{0}_{\ell^*(A)}(S^1;su(3))$ (note that its left translation to
$T_{H_{\ell}}SU(3)$ is tangent to the adjoint orbit of $H_{\ell}$).
\end{pf}

The following `warm-up' proposition treats the case $SU(2).$

\begin{prop} \label{6.1}
If $A$ is an irreducible flat $SU(2)$ connection and $a \in \HH^1_A(X; su(2))$
is nonzero,
then there exists a curve $\ga$ with $\tr(H_\ga I_\ga) \neq 0.$
\end{prop}
\begin{pf}
Since $a$ is nonzero and harmonic, by Proposition \ref{first loop}
we can choose a curve $\ell$ so that $\Pi_{z(H_{\ell})}(
I_\ell) \neq 0.$
Gauge transform $A$ so that
$$H_\ell =  \left( \begin{array}{cc} \la & 0 \\ 0 & \la^{-1}
\end{array}\right)$$
is diagonal and write
$$I_\ell = \left( \begin{array}{cc} i \al & \be \\ \bar{\be} & -i \al
\end{array} \right), \quad 0 \neq \al \in \RR.$$
Then $\tr(H_\ell I_\ell)= \al(\bar{\la} - \la) \neq 0$
unless $H_\ell = \pm I.$
Taking $\ga = \ell$ proves the claim if $H_\ell \neq \pm I.$
Otherwise, we can always find $\ga$ so that $\tr(H_\ga I_\ell) \neq 0.$
Using the fact that $H_\ell$ is central,  it follows that
\begin{eqnarray*}
\tr(H_{\ga \cdot \ell} I_{\ga \cdot \ell})
&=&\tr(H_\ga H_\ell I_\ell) +\tr(H_\ga I_\ga H_\ell) \\
&=& \pm (\tr(H_\ga I_\ell) +\tr(H_\ga I_\ga)).
\end{eqnarray*}
Since $\tr(H_\ga I_\ell)$ is nonzero, it follows that either
$\tr(H_{\ga \cdot \ell} I_{\ga \cdot \ell})$ or
$\tr(H_\ga I_\ga)$ is also nonzero, and this proves the proposition.
\end{pf}

The  same is true for $SU(3),$ but it takes more work to prove.  
\begin{thm} \label{6.2}
If $A$ is an irreducible flat $SU(3)$
connection and $a \in \HH^1 _{A}(X; su(3))$ is nonzero,
then there exists a curve $\ga$ with $\tr(H_\ga I_\ga) \neq 0.$
\end{thm}
\begin{pf}
Choose $\ell$ so that $\Pi_{z(H_{\ell})}(I_\ell) \neq 0.$
Gauge transform  so that
$$H_\ell = \left( \begin{array}{ccc} \la_1 && 0 \\  &\la_2 \\ 0&& \la_3
\end{array}\right)$$
is diagonal and write
$$I_\ell = \left( \begin{array}{ccc} i \al_1 &&  * \\ & i\al_2 \\
\bar{*} && i\al_3  \end{array} \right),$$
where $\al_{i}$ are real numbers, not all zero.
Of course, $\la_3 = (\la_1 \la_2)^{-1}$ and $\al_3 = -\al_1 -\al_2.$

If $H_\ell$ has only one eigenvalue, namely if $\la_1=\la_2=\la_3,$
then $H_\ell$ is central and the theorem follows from the same argument as
was used to
prove Proposition \ref{6.1}.
Otherwise, either $H_\ell$ has three distinct eigenvalues or it may be
further conjugated
so that $\la_1=\la_2 \neq \la_3.$
The following argument treats only the first of these two cases.
The second case
requires a more elaborate argument, given in the next subsection.

Assume $\la_1, \la_2$ and $\la_3$ are all distinct.
Suppose first of all that $\al_i = 0$ for some $i,$ which can be taken (wlog)
to be $i=3.$ Since $\tr(I_\ell)=0,$
$$\tr(H_\ell I_\ell) = \la_1 (i\al_1) +\la_2(i\al_2) = i \al_1 ( \la_1
-\la_2),$$
which is nonzero since $\al_1 \neq 0$ and $\la_1 \neq \la_2.$

Now suppose $\al_i \neq 0$ for all $i.$ By replacing $a$ with $-a,$ if
necessary, we can assume that
two of the $\al_i$'s are positive, which we take (wlog) to be $\al_1$ and
$\al_2.$
Then
$$\tr(H_\ell I_\ell) = i \la_1 \al_1 + i \la_2 \al_2 - i(\la_1
\la_2)^{-1}(\al_1+\al_2).$$
Thus $\tr(H_\ell I_\ell) =0$ implies
$\la_1 \al_1 +  \la_2 \al_2 = (\la_1 \la_2)^{-1}(\al_1+\al_2).$
If this were the case, then $| \la_1 \al_1 +  \la_2 \al_2| = |\al_1 + \al_2|,$
which is only possible if $\la_1 = \la_2,$ a contradiction.

The following example illustrates the difficulties
when $H_\ell$ has only two distinct eigenvalues.
Suppose
$$H_\ell = \left( \begin{array}{ccc} \la &&0\\ & \la \\ 0 && \la^{-2}
\end{array} \right).$$
Then one see that $\tr(H_\ell I_\ell) =0$ if
$$I_\ell = \left(\begin{array}{ccc} i \al && * \\ & -i \al \\ \bar{*} && 0
\end{array} \right).$$

The next subsection is devoted to treating this problematic case.
Observe that we can assume that $H_\ell$ has infinite order
for the following reason.
If $H_\ell$ has finite order $k$
and if $\ga$ is chosen 
so that $\tr(H_\ga I_{\ell^k}) \neq 0,$
then just as in the proof of Proposition \ref{6.1}, we compute that
$$\tr(H_{\ga \cdot \ell^k} I_{\ga \cdot \ell^k}) =\tr(H_\ga I_{\ell^k}) +
\tr(H_\ga I_\ga).$$
But $\tr(H_\ga I_{\ell^k}) \neq 0,$ hence it follows that one of the other two
terms is also non-zero.

\subsection{Linear algebra} \label{linear algebra}
In this subsection, we
complete the proof of Theorem \ref{6.2} demonstrating
the existence of perturbation curves with certain properties.
The remaining case is when $H_\ell$ 
has only two distinct eigenvalues.
As indicated in the previous subsection,
we can further assume that
$H_\ell$  has infinite order.
Although $H_\ell$ may not have three distinct eigenvalues,
the following proposition assures us that $H_\ga$ has three
distinct eigenvalues for some loop $\ga.$

\begin{prop} \label{6.3}
If $\varrho: \pi_1(X) \lto SU(3)$ is an irreducible representation,
then there exists some element $\ga \in \pi_1(X)$ such that $\varrho(\ga)$
has three distinct eigenvalues.
\end{prop}
\noindent
{\it Remark.} Besides the existence of the irreducible, rank three
representation,
the proof makes no assumptions on the group $\pi_1(X).$
\begin{pf}
By irreducibility, we can find $\ell$ with $\varrho(\ell)$ noncentral.
Set $L = \varrho(\ell).$
Obviously, we are done unless $L$ has only two distinct eigenvalues.
Since the conclusion of the proposition is invariant under conjugation,
we can assume
$$L = \left( \begin{array}{ccc} \la &&0\\ & \la \\ 0 && \la^{-2}
\end{array} \right)$$
is diagonal.
By irreducibility of $\varrho,$ there exists $m \in \pi_1(X)$ so that
$\varrho(m)$ does not commute with $L.$ Set $M = \varrho(m).$
Thus neither $M$ nor $LM$ is diagonal.
Of course, we can also assume that $M$ has only two distinct eigenvalues;
otherwise we are done! 
 Let $\mu$ be the eigenvalue of $M$ of multiplicity two.
Now both $L$ and $M$ have 2-dimensional eigenspaces, so for dimensional
reasons,
$L$ and $M$ have a common eigenvector.
After conjugating by a matrix commuting
with $L,$  it follows that $M$ can be written in block diagonal form:
$$M=\left( \begin{array}{cc} \mu & 0  \\ 0 & A \end{array} \right)$$
where
$$A=\left( \begin{array}{rr}a & b \\ -\bar{b}& \bar{a}  \end{array}\right)
\left( \begin{array}{cc}\mu & 0 \\ 0 & \mu^{-2} \end{array}\right)
\left( \begin{array}{rr}\bar{a} & -b \\ \bar{b}&  a  \end{array}\right)$$
and $|a|^2+|b|^2=1.$
The matrix product $LM$ also comes in block form:
$$LM=\left( \begin{array}{cc} \la \mu  & 0 \\ 0 & B \end{array} \right)
\quad \hbox{ where } \quad
B = \left( \begin{array}{cc}\la & 0\\ 0 & \la^{-2}\end{array}\right) A.$$
We claim that $LM$ has three distinct eigenvalues.
First of all, notice that the two eigenvalues of $B$ are distinct; 
otherwise $LM$ would be diagonal, in which case $L$ and $M$ would commute.
So it suffices to prove that $\la \mu$ is not an eigenvalue of $B.$

Suppose to the contrary that $\la \mu$ is an eigenvalue of $B$,
i.e., suppose
$(\la\mu)^2 -\tr(B) (\la \mu) + \det(B) =0.$
Computing $\tr(B)$ directly,
one finds
$$\tr(B)  = |a|^2(\la \mu + \la^{-2} \mu^{-2}) + (1-|a|^2)(\la \mu^{-2} +
\la^{-2} \mu).$$
Plugging this into the characteristic equation   and using $\det B =
\la^{-1} \mu^{-1}$
gives
$$(1-|a|^2)(\la^2 \mu^2 + \la^{-1}\mu^{-1} - \la^2\mu^{-1}-\la^{-1}\mu^2) =
0.$$
So either $|a|= 1,$  implying $A=\pm I$ and contradicting our choice of $M,$
or, after multiplying by $\la \mu,$
\begin{eqnarray*}
0 &=& \la^3 \mu^3 + 1 - \la^3-\mu^3
= (\la^3 -1)(\mu^3 - 1).
\end{eqnarray*}
However,
$\la^3 =1$ implies $L$ is central and  $\mu^3 = 1$ implies
$M$ is central, each giving contradictions. Hence
$\la \mu$ is not an eigenvalue of $B,$ which proves our claim.
\end{pf}

With regard to Theorem \ref{6.2},
we have already proved the existence of $\ga$ unless $H_\ell(A)$ has two
distinct eigenvalues, so assume
$$H_\ell(A) = \left( \begin{array}{ccc} \la && 0 \\ & \la \\ 0 && \la^{-2}
\end{array} \right).$$
Set $L_t = H_\ell(A_t)$ where $A_t = A + t a.$
Observe moreover that we are done unless
$$\left. \frac{d L_t}{dt}\right|_{t=0} = L_0 \left( \begin{array}{ccc} i \al
&&0\\ & -i \al\\ 0 && 0 \end{array} \right),$$
where $\al \neq 0.$

Before stating the next lemma, we make a definition.

\begin{defn} \label{6.4}
For a fixed angle $\eta \in [0, 2 \pi),$ let $G_\eta$ be the subset of $SU(3)$
consisting of matrices of the form
$$M= \left( \begin{array}{lll}
			a &  		b e^{-i \eta} & 	c \\
			be^{i \eta} & a & 			c e^{i \eta} \\
			c'&		c' e^{-i \eta} &	d
\end{array} \right)$$
for $a,b,c,c',d \in {\Bbb C}.$ Note that $M \in SU(3) \Rightarrow |c|=|c'|.$
\end{defn}
Clearly $G_\eta$ is a subgroup; in fact,
matrices of this form in $SL(3,\CC)$ form a subgroup of complex codimension 4,
so one would expect that $G_\eta$ has real codimension 4 in $SU(3).$
This is indeed the case; if $M$ is chosen as above and
$$P=\frac{1}{\sqrt{2}}\left( \begin{array}{ccc}
			1 &  		0 & 	1 \\
			e^{i \eta} & 0 & 	-e^{i \eta} \\
			0&		\sqrt{2} &	0   \end{array}
\right),$$
then $P \in U(3),$  
$$P^{-1}M P = \left( \begin{array}{ccc}
			a+b &  		\sqrt{2} c & 	0 \\
			\sqrt{2} c' & d & 			0 \\
			0&		0 &	a-b   \end{array} \right),$$
and so $G_\eta$ is conjugate to the subgroup $S(U(2)\times U(1)).$

\begin{lem} \label{6.5}
Suppose that $L_t, M_t:(-\ep,\ep) \lto SU(3)$  are smooth paths. 
Write $L'_0 = L_0^{-1} \left. \frac{dL_t}{dt} \right|_{t=0},$ and assume both 
$L_0$ and $L'_0$ are
diagonal, with
$$L_0  =\left( \begin{array}{ccc} \la  && 0 \\ & \la  \\ 0 && \la^{-2}
\end{array} \right) \hbox{ and }
L'_0 =  
\left( \begin{array}{ccc} i\al  & & 0 \\ & -i\al  \\ 0 & & 0
\end{array} \right)$$
for $\la$ a complex unit of infinite
order and for $\al \neq 0$.
If
$$\left. {\tfrac{d}{dt}} \tr(W_t) \right|_{t=0} = 0$$
for every word $W_t$ in $L_t$ and $M_t$,
then $M_0 \in G_\eta$ for some $\eta.$

\end{lem}

\begin{pf}
For general $L,M \in SU(3)$ with
$$L = \left( \begin{array}{ccc} \la_1  &&0\\ & \la_2 & \\ 0 && \la_3
\end{array} \right)$$
 diagonal and
$M = (\mu_{ij})$  arbitrary, it is not difficult to verify that
\begin{eqnarray*}
\tr(LM) &=& \sum_{i=1}^3 \la_i  \mu_{ii} , \\
\tr(LML^{-1} M^{-1}) &=& \sum_{i,j=1}^3 \la_i  \bar{\la}_j |\mu_{ij}|^2.
\end{eqnarray*}

Now suppose $L_t, M_t$ are as in the hypotheses. We write $M_t = (\mu_{ij}(t))$  
and let $\mu_{ij} = \mu_{ij}(0)$ for
convenience.
Applying the above formula to $L_t^k M_t$ and $L_t^k M_t L_t^{-k} M_t^{-1}$
and taking derivatives,
we see from the hypotheses that
\begin{eqnarray*}
0 &=& 
\left. {\tfrac{d}{dt}} \tr(L_t^k M_t) \right|_{t=0} \\
&=& \la^k  \left. {\tfrac{d}{dt}} 
\left(\mu_{11}(t) + \mu_{22}(t)\right)\right|_{t=0} +
\bar{\la}^{2k}  \left. {\tfrac{d}{dt}}\mu_{33}(t)\right|_{t=0}
 + i k \al \la^k \left(\mu_{11}- \mu_{22}\right),
\end{eqnarray*}
and that
\begin{eqnarray*}
0 &=& \left. {\tfrac{d}{dt}} \tr(L_t^k M_t L_t^{-k} M_t^{-1}) \right|_{t=0} \\
&=& \left. {\tfrac{d}{dt}} \left(|\mu_{11}(t)|^2 + |\mu_{12}(t)|^2+ |\mu_{21}(t)|^2+
|\mu_{22}(t)|^2+ |\mu_{33}(t)|^2\right) \right|_{t=0}  \\
& & + \; \la^{3k} \left. {\tfrac{d}{dt}} \left(|\mu_{13}(t)|^2 + |\mu_{23}(t)|^2\right) \right|_{t=0} +
 \bar{\la}^{3k} \left. {\tfrac{d}{dt}} \left(|\mu_{31}(t)|^2 + |\mu_{32}(t)|^2\right) \right|_{t=0} \\
& & + \; k \al\left\{ 2(|\mu_{12}|^2 - |\mu_{21}|^2)
 + \la^{3k} (|\mu_{13}|^2 - |\mu_{23}|^2)
-\bar{\la}^{3k} (|\mu_{31}|^2 - |\mu_{32}|^2) \right\}.
\end{eqnarray*}
Since both equations hold for all $k \geq 0$ and since $\la$ has infinite
order, we deduce that:
$$\begin{array}{ll}
\; {\rm (i)} \; \; \mu_{11} = \mu_{22},   \quad \quad &  {\rm (ii)} \; \; |\mu_{12}| = |\mu_{21}|,\\
{\rm (iii)} \; \; |\mu_{13}| = |\mu_{23}|, \quad \quad   & {\rm (iv)} \; \; |\mu_{31}| = |\mu_{32}|.
\end{array}$$
Here, (i) is a consequence of the first equation and  (ii)--(iv) come from
the second.
The last three conditions are equivalent to the existence of angles
$\eta_1,\eta_2,$ and $\eta_3$ with
 $\mu_{21} = e^{i 2 \eta_1} \mu_{12}, \;\;  \mu_{23} = e^{i \eta_2} \mu_{13},\;$
and $\mu_{32}= e^{-i \eta_3} \mu_{31}.$
 
To conclude that $M_0 \in G_\eta,$ we just need to
show that $\eta_1 = \eta_2 = \eta_3 \mod (2\pi).$
Applying (i) to $(M_0)^2$ implies
$\mu_{13} \mu_{31} = \mu_{23} \mu_{32},$ thus
$\eta_2 = \eta_3 \mod(2 \pi).$
Now apply the unitary condition to $M_0$ to see
$0=\sum_{j=1}^3 \mu_{ij} \bar{\mu}_{3j}$ for $i = 1,2.$ Comparing
these, we conclude $\eta_1 = \eta_2 \mod(2 \pi).$
This completes the proof of the lemma.
\end{pf}

To establish Theorem \ref{6.2}, we seek a curve
$\ga$ such that $\tr (H_\ga I_\ga) \neq 0.$
Setting $A_t = A+ta,$ this is equivalent to the condition that
$\left. \frac{d}{dt} \tr H_\ga(A_t) \right|_{t=0} \neq 0.$
According to the previous lemma,
letting $\ga$ range over all words in $L_0$ and $M_0,$ the 
only way this can fail 
is if $M_0 \in G_\eta$ for some $\eta.$
We shall show in the following argument that
the irreducibility of
$A$ guarantees the existence of an $M = H_{m}(A)$ such that
$M \not\in G_\eta$ for any $\eta.$

{\it Proof of Theorem \ref{6.2}.} \;
We provide the proof in the remaining case when
$L_0 = H_\ell(A)$ has two distinct eigenvalues and is of infinite order.
Set $A_t = A+ta$ and $L_t = H_\ell(A_t).$
By Proposition \ref{6.3}, we have a loop $m_1$
such that
$H_{m_1}(A) $ has three distinct eigenvalues. 
Set $M_1 = H_{m_1}(A)$ and $M_{1,t} = H_{m_1}(A_t)$.
Assume first that $M_1 \not\in G_\eta$ for any $\eta \in [0,2 \pi].$
By Lemma \ref{6.5}, there is a word $W_t$ in $L_t$ and $M_{1,t}$
such that $  \left. \frac{d}{dt} \tr W_t \right|_{t=0}  \neq 0.$
Taking $\ga$ as the loop obtained from the corresponding word in $\ell$
and $m_1$, then $W_t = H_\ga(A_t)$ and hence
$\left. \frac{d}{dt} \tr H_\ga(A_t) \right|_{t=0}  \neq 0,$
which proves the theorem in this case.

So now suppose $M_1 \in G_{\eta_1}$ and write
$$M_1= \left( \begin{array}{lll}
			a_1 &  		b_1 e^{-i \eta_1} & 	c_1 \\
			b_1e^{i \eta_1} & a_1 & 		c_1 e^{i
\eta_1} \\
			c'_1&		c'_1 e^{-i \eta_1} &	d_1
\end{array} \right).$$
Although $A_t$ has been gauge transformed
so that the path $H_\ell(A_t)$ is diagonal,
we can further conjugate by a diagonal matrix since it
acts trivially on $H_\ell(A_t)$.
Applying such a conjugation to $M_1,$ we can arrange 
that $b_1$ and $c_1$ are both real and non-negative.

Since $A$ is irreducible, we can choose $m_2 \in \pi_1(X)$
such that  $M_2=H_{m_2}(A) \not\in G_{\eta_1}.$
Repeating the argument above with $M_1$ replaced by $M_2,$
we can assume that $M_2 \in G_{\eta_2}$ for some
$\eta_2$ and write
$$M_2= \left( \begin{array}{lll}
			a_2 &  		b_2 e^{-i \eta_2} & 	c_2 \\
			b_2e^{i \eta_2} & a_2 & 		c_2 e^{i
\eta_2} \\
			c'_2&		c'_2 e^{-i \eta_2} &	d_2
\end{array} \right).$$

We now claim that one of $M_1 M_2$ and $M_1^{-1} M_2$ is not contained in
$G_\eta$ for any $\eta.$
This completes the proof of Theorem \ref{6.2} by 
repeating once
again the above argument 
with $M_1$ replaced by either $M_1 M_2$ or $M_1^{-1} M_2$
and invoking Lemma \ref{6.5} to produce the curve $\ga$ with the
desired properties.

So, it only remains to establish the claim, which is 
proved by contradiction.
Suppose that $M_1 M_2$ is contained in $G_\eta$ for some $\eta.$
Equating the $(1,1)$ and $(2,2)$ entries of $M_1 M_2,$ we find
\begin{equation} \label{e6.1}
b_1 b_2 u^2 + c_1 c_2' u  - (b_1 b_2 + c_1 c_2')=0,
\end{equation}
where $u = e^{i(\eta_2-\eta_1)}.$

Suppose first of all that $b_1 = 0.$ Because $c_1 = 0 \Rightarrow M_1$ is
diagonal,
and $u=1 \Rightarrow M_2 \in G_{\eta_1},$ neither of which is the case,
the only possibility is that  $c_2'=0.$
But writing out $M_1 M_2$
and demanding that  the off-diagonal
terms have the required form, this would imply that $b_2 =0$ or $u=1,$
both of which lead to contradictions.

So assume $b_1 \neq 0.$
By similar considerations, we can also assume $b_2 \neq 0.$
Solving equation (\ref{e6.1}) for $u$ gives
$$u = -1 - \frac{c_1 c_2'}{b_1 b_2},$$
since we have already seen that the other possibility,
namely $u=1,$ leads to a contradiction.

The same reasoning applied to $M_1^{-1} M_2$ shows that
$$u = -1 -\frac{\bar{c}'_1 c_2'}{\overline{b}_1 b_2}.$$
Equating these two formulas for $u$ gives
$$c_1 \overline{b}_1 = \overline{c}'_1 b_1.$$
Since $b_1$ and $c_1$ are real, this shows that
$\overline{c}_1' = c_1 = c_1',$  which forces both
$a_1$ and $d_1$ to also be real.
It immediately follows that $M_1^{-1} = M_1^* = M_1,$
hence the eigenvalues of $M_1$ equal $\pm 1.$ In particular,
$M_1$ has at most two distinct eigenvalues, which
contradicts our choice
of $M_1.$
This proves the claim and concludes the proof of
Theorem \ref{6.2}.
\end{pf}
Given loops $\ell_1,\ldots, \ell_n$ in $X,$ define gauge invariant functions
$f_j, g_j:\AA \lto \RR$ for $j=1,\ldots, n$ to be the
real and imaginary parts of $\tr \hol_{\ell_j}(A),$ so that
$$\tr(\hol_{\ell_j}(A)) = f_j(A) + i g_j(A).$$
Note that if $A$ is an $SU(2)$ connection, then $g_j(A)=0.$

\begin{cor}\label{florida} 
\begin{enumerate}
\item[(i)] If $A$ is an irreducible, flat $SU(2)$ connection, then
there exist loops $\ell_1,\ldots, \ell_n$
so that the map 
from $\HH^1_A(X;su(2))$ to $\RR^n$ 
given by \\
$a \mapsto (D f_1 (A)(a),\ldots, D f_n (A)(a))$
is injective. 
\item[(ii)] If $A$ is an 
irreducible, flat $SU(3)$ connection, then there exist loops
$\ell_1,\ldots, \ell_n$
so that the map
from $\HH^1_A(X;su(3))$ to $\RR^{2n}$
given by \\
$a \mapsto (D f_1 (A)(a), D g_1 (A) (a), \ldots , D f_n (A)(a), D g_n(A)(a))$
is injective.
\end{enumerate}
\end{cor}

\subsection{Second order arguments} 
\label{s6.3}

Suppose now that $A$ is a reducible flat $SU(3)$ connection on $X$.
Then part (i) of Corollary \ref{florida} allows us to find loops about which the
real part of the
derivative of trace of holonomy detects any first order deformations
of $A$ in directions tangent to the reducible stratum,
i.e., in the directions of $\HH ^{1 }_{A}(X; \rr)$. But
invariance  under the gauge group,  in particular under
$\stab (A) = U(1)$, prevents the derivative from detecting first order deformations
in directions normal to the reducible stratum,
i.e., in the directions of $\HH ^{1 }_{A}(X;\rrp).$    
Instead, we consider
second derivatives of the gauge invariant functions in these
directions.  This portion of the argument  closely parallels
the argument used to handle abelian flat connections in
the $SU(2)$ moduli space \cite{herald1}.  

Notice first that $\HH^1_A(X;\rrp)$ is
a module over the quaternions $\BH.$
To see this, let $SP(1)$ be the unit quaternions and define 
$\phi: SU(2) \rightarrow SP(1)$
by
$$\left(
\begin{array}{rr}
 a & -\bar{b} \\
b & \bar{a}
\end{array} \right) \mapsto a+Jb$$
and $F:\CC^2 \to \BH$ by
$F(v_1,v_2) = v_1+Jv_2.$
Then for $A \in SU(2)$ and $v \in \CC^2,$
$$F(Av) = \phi(A) F(v).$$
This turns the action of $SU(2)$ on $\CC^2$ into
left multiplication by elements of $SP(1)$ on $\BH.$

Now suppose
$\varrho: \pi_1(X) \to
SU(2)$ is
an irreducible representation and 
let
$E_\varrho$ be the flat bundle $\widetilde{X} \times_{\pi_1(X)} \CC^2$,
where
$\widetilde{X}$ is the universal cover of $X$ and
$\pi_1(X)$ acts by deck transformations on $\widetilde{X}$
and via the canonical representation of $\varrho$ on $\CC^2.$
We identify $E_\varrho$ as a flat bundle with the
subbundle of $\ad P = X \times su(3)$ corresponding to $\rrp \subset su(3).$
The de Rham theorem provides an isomorphism
$\HH^1_A(X;\rrp) \cong H^1(X;E_\varrho).$
Here, $H^1(X;E_\varrho) = Z^1(X;E_\varrho)/ B^1(X;E_\varrho)$
is by definition the space of 1-cocycles
modulo the 1-coboundaries.
Using a presentation
$\pi_1(X) = \langle x_1,\ldots, x_n \mid r_1,\ldots,r_m
\rangle,$
we can identify
the 1-cochains as elements
$(v_1,\ldots, v_n) \in \CC^2 \times \cdots \times \CC^2
\cong \BH^n$
and the subspaces 
$Z^1(X;E_\varrho)$ of 1-cocycles and $B^1(X;E_\varrho)$ of
1-coboundaries as
submodules.
For example, $(v_1,\ldots, v_n)$ is a
coboundary
if and only if there is some $v \in \CC^2$
such that
$v_i = v - \varrho(x_i) v$ for $i=1,\ldots,n.$
Observe that $B^1(X;E_\varrho)$ is closed under
right multiplication by elements in $\BH.$
Similarly, $(v_1,\ldots, v_n)$ is a cocycle
if and only if
the following linear equations, which are derived from the relations
$r_1,\ldots, r_m$ using the Fox differential calculus,
are satisfied:
\begin{eqnarray} \label{star}
&M_{11} v_1 + \cdots + M_{1n} v_n = 0& \nonumber \\
& \vdots & \\
&M_{m1} v_1 + \cdots + M_{mn} v_n = 0.& \nonumber
\end{eqnarray}
Here each $M_{ij}$ is a sum of $SU(2)$ matrices and thus is a
$2 \times 2$ matrix of the form
$$M_{ij} = \left( \begin{array}{cc} a_{ij}  & - \bar{b}_{ij} \\
b_{ij} & \bar{a}_{ij}
\end{array} \right)$$
for some $a_{ij}, b_{ij} \in \CC.$
For $v=(v_1,v_2)
\in \CC^2$ and
$h = z_1+J z_2 \in \BH,$ where $z_1,z_2 \in \CC,$
set
$v \cdot h = ( z_1 v_1 - z_2 \bar{v}_2, z_2 \bar{v}_1 + z_1 v_2) \in \CC^2$.
(This is just multiplication in $\BH$
under the isomorphism $F:\CC^2 \cong \BH.$)
Now if $(v_1,\ldots, v_n)$ satisfies (\ref{star}) above, then so does
$(v_1 \cdot h, \ldots, v_n \cdot h).$
This shows that $Z^1(X;E_\varrho)$ is closed under
right multiplication by elements of $\BH.$
Since both $B^1(X;E_\varrho)$ and $Z^1(X;E_\varrho)$
are right $\BH$-modules, so is
$H^1(X;E_\varrho) = Z^1(X;E_\varrho)/B^1(X;E_\varrho).$

Recalling the notation of
Definition \ref{hermite}, we use $\sym \HH^{1}_{A}(X;\rrp)$
to  denote the set
of all symmetric $\stab (A) \cong U(1)$ invariant bilinear forms on
$\HH^{1 }_{A}(X;\rrp)$, regarded as a real vector space with a
$U(1)$ action. 

\begin{prop}  \label{hessian loops}  If $A$ is a reducible flat connection,
then there exist loops $\ell_1, \ldots , \ell_{n_1}$ in $X$ and a set
$F = \{ f_1,\ldots,f_n\}$ of gauge invariant functions such that:
\begin{enumerate}
\item[(i)]
Each $f_i \in F$ is the real or imaginary part of $\tr \hol_{\ell_j}$
for some $j=1,\ldots,n_1.$
\item[(ii)] 
The map
$\RR^{n} \lto \sym \HH^{1}_{A}(X;\rrp)$ given by 
$ (x_1, \ldots, x_{n})\mapsto  
\sum_{i=1}^n x_i \hess f_i(A) $
is surjective.
\item[(iii)] $Df_i(A) =0$ for $i=1,\ldots,n.$
\end{enumerate}
\end{prop}
\begin{pf}
Assume $A$ has been gauge transformed to
take values in $su(2) \subset su(3)$ and denote by $\widehat{A}$
the associated irreducible $SU(2)$ connection.
In order to construct the loops $\ell_{1}, \ldots , \ell_{n_1}$,
we will need to introduce curves in $X$ that are in a certain sense
dual to a basis for $\HH^{1}_A(X;\rrp)$ over $\BH$.

Let
$\varrho:\pi_1(X) \to SU(2)$  be the irreducible
$SU(2)$ representation associated to $\widehat{A}$,
and let 
$E_\varrho = \widetilde{X} \times_{\pi_1(X)} \CC^2$
as before.
Consider $H_i(X; E_\varrho)$, 
homology with local coefficients in $E_\varrho,$
which is by definition the homology of the complex 
$$\cdots \lto C_i(\widetilde{X}) \otimes_{{\Bbb Z}[\pi_1(X)]} \CC^2 
\stackrel{\partial_i \otimes 1}{\lto} C_{i-1}(\widetilde{X}) 
\otimes_{{\Bbb Z}[\pi_1(X)]} \CC^2 \lto \cdots.$$
From our previous discussion,
it is not hard to see that $H_1(X;E_\varrho)$ is
a right $\BH$-module.
Thus, we have a basis for
$H_{1}(X; E_\varrho)$ over $\BH$
consisting of classes each of which can
be represented by a $\CC^2$-labelled curve $\widetilde{\ga}_i$ in
the universal cover $\widetilde{X}$ of $X.$
Each $\widetilde{\ga}_i$ is a lift of a loop $\ga_i$ in $X$
with
$\hol_{\ga_i}(A) =1$
(because
the labelled lift of $\ga_i$
lies in 
$\ker \partial_1 \otimes 1 $ if and only if
the holonomy of $A$ around $\ga_i$ is trivial).
Let $\om_{1}, \ldots , \om_{m}$ be the Hom dual basis for
$H^1(X;E_\varrho)$ over $\BH$.
Of course, each $\om_i$ determines a real, 4-dimensional
subspace $V_i = \phi( \Span_\BH \om_i) \subset \HH_A^1(X;\rrp)$
where $\phi: H^1(X;E_\varrho) \to \HH_A^1(X;\rrp)$ 
is the isomorphism provided by the de Rham theorem.
Each $V_i$ is preserved by
the subgroup $\stab(A) \subset \GG$, thus 
\begin{equation} \label{decomp}
\HH_A^1(X;\rrp) = V_1 \oplus \cdots \oplus V_m
\end{equation}
is a decomposition of $\HH^1_A(X;\rrp)$
into 2-dimensional complex vector spaces.
We denote by $a_i$ the image of $a \in \HH_A^1(X;\rrp)$ 
under the projection $p_i : \HH_A^1(X;\rrp) \to V_i.$

Let $\UU = \sym \HH^1_A(X;\rrp)$ be
the space of symmetric $\stab(A) = U(1)$ invariant bilinear forms
on $\HH^1_A (X; \rrp)$.
Our goal is to find a collection of loops such that
the Hessians of the real and imaginary parts of the trace
of holonomy functions around these loops span $\UU$.

There is a decomposition
of $\UU$ corresponding
to (\ref{decomp})
given by 
$\UU =\bigoplus_{i \leq j} \; \UU_{ij}$,
where $B \in \UU_{ij}$ in case
$$B(a,b) = \left\{ 
\begin{array}{ll} B(a_i,b_i)  & \hbox{ if } i=j, \\
B(a_i, b_j)+B(a_j,b_i) & \hbox{ if } i \neq j.
\end{array} \right.$$
Thus  every
$ B \in \UU_{ij}$ is entirely determined by its restriction to
$V_i \times V_j$. Let $\{ a, b\}$ be a basis for $V_i$
and $\{ c,d\}$ a basis for $V_j$. 
In terms of the real bases 
$\{ a, ia, b, ib \}$ for $V_i$ and $\{ c, ic, d, id \}$
for $V_j$, the restriction of 
$B$ to $V_i \times V_j$ is a real 
$4 \times 4$ matrix of the form
\begin{equation} \label{Uform}
\left[ \begin{array}{rrrr} 
x & 0 & y & z \\
0 & x & -z & y \\
y & -z & w & 0 \\
z & z & 0 & w \\
\end{array} \right]
\quad \hbox{ if $i=j,$ and}
\quad \quad
\left[ \begin{array}{rrrr}
p & q & r & s \\
-q & p & -s & r \\
t & -u & v & w \\
u & t & -w & v \\
\end{array} \right] \hbox{ if $i \neq j$.}
\end{equation}
From this, it follows that
$$\dim \UU_{ij} = \left\{
\begin{array}{ll} 4 & \hbox{ if } i=j \\ 8 & \hbox{ if } i \neq j.
\end{array} \right.$$
We prove the proposition by 
constructing, for each $i \leq j,$
gauge invariant functions 
satisfying conditions (i) and (iii) such that
their
Hessians at $A$ span
$\UU_{ij}$.
We begin with the case $i = j$.

Given $\be:[0,1] \to X$ with $\be(0)=\be(1),$
parallel translation can be used to associate a function $\al:[0,1] \to \rrp$
to any $su(3)$-valued 1-form $a$ by setting 
\begin{equation} \label{funk}
\al(t)dt = P_\be(0,t)^{-1} a_{\be(t)} P_\be(0,t),
\end{equation}
where $P_\be(0,t)$ is parallel translation by $A$ along
$\be$ from $\be(0)$ to $\be(t)$.
If $a \in \HH^1_A(X; \rrp),$  then
$\int_0^1 \al(t) dt =\int_{\ga_i} a \in \rrp $.
The linear transformation
$\HH^1_A(X; \rrp) \to \rrp$
defined by $a \mapsto \int_{\ga_i} a$ 
has kernel $V_1 \oplus \cdots \widehat{V}_i \cdots \oplus V_m$
(because the basis
$\om_1,\ldots, \om_m$ is Hom dual to 
$\widetilde{\ga}_1, \ldots, \widetilde{\ga}_m$) and
determines an isomorphism
$V_i \to \rrp.$

Note how the correspondence (\ref{funk}) behaves for products of loops.
If $\be = \ell_1 \cdots \ell_k : [0,k]\to X$
(where each $\ell_i:[i-1,i] \to X$ is a
loop),
define $\al_i:[i-1,i] \to \rrp$ by
$\al_i(t) dt = P_{\ell_i}(i-1,t)^{-1} a_{\ell_i(t)} P_{\ell_i}(i-1,t).$
Defining $\al:[0,k] \to \rrp$ by (\ref{funk}), then
\begin{equation} \label{funky}
\al(t)   = 
hol_{\ell_{1}}(A)^{-1} \cdots hol_{\ell_{i-1}}(A)^{-1} \al_i(t) 
\; hol_{\ell_{i-1}}(A) \cdots hol_{\ell_{1}}(A)
\end{equation}
for $t \in [i-1,i].$
\begin{lem} \label{loopy-0} 
Suppose $\ell$ is a loop with $L:=\hol_\ell(A) \in SU(3)$ nontrivial.
If
$a,b \in \HH^1_A(X;\rrp)$
and if we
set $\xi_i = \int_{\ga_{i}} a  \in \rrp$ 
and $\zeta_i = \int_{\ga_{i}} b \in \rrp,$
then
\begin{itemize}
\item [(i)]  $\hess \tr \hol_{\ga_i}(A)(a,b)= 2 \tr(\xi_i \zeta_i).$
\item [(ii)]  $\hess \tr \hol_{\ell \cdot \ga_i}(A)(a,b)= 
\tr(L (\xi_i \zeta_i + \zeta_i \xi_i)).$
\end{itemize}
\end{lem}

\begin{pf}
Let $B, B_\ell:\HH^1_A(X;\rrp) \times \HH^1_A(X;\rrp) \to \CC$
be the symmetric, bilinear pairings
coming from the Hessians at $A$ of 
$\tr \hol_{\ga_i}$ and $\tr \hol_{\ell \cdot \ga_i},$ respectively.
(Notice that $B(a,b) \in \RR$ 
because $\hol_{\ga_i}(A)$ is trivial. This
follows from \ref{derivative of tr hol} (ii) 
and the elementary fact that
$\tr(\xi \zeta ) \in \RR$ for $\xi, \zeta \in \rrp$.)

Since $B$ and $B_\ell$ are both symmetric, and since
$\int_{\ga_i} (a+b) = \xi_i + \zeta_i,$
it suffices to show (i) and (ii)
in the case $a=b.$
To prove (i),
parameterize $\ga_i$ by the interval $[0,1]$ 
and define $\al: [0,1] \to \rrp$ 
associated to the 1-form $a \in \HH^1_A(X; \rrp)$
using (\ref{funk}). 
Using the formula from  
Corollary \ref{derivative of tr hol} (ii)
and noting that $\hol_{\ga_i}(A)$ is trivial,
it follows that
\begin{eqnarray*} 
\hess \tr \hol_{\ga_i}(A)(a,a) 
&=& \int_0^1 \int_0^s \tr (\al(s) \al(t) + \al(t) \al(s))  dt ds\\
&=& \int_0^1 \int_0^1 \tr (\al(s) \al(t))  dt ds\\
&=& \tr \left(\int_0^1 \al(s) ds \int_0^1 \al(t) dt \right) 
= \tr (\xi_i^2).
\end{eqnarray*}
This proves (i).

To prove (ii), set $\be = \ell \cdot \ga_i$ and parameterize it
by the interval $[0,2]$ so that the subintervals
$[0,1]$ and $[1,2]$ parameterize $\ell$ and 
$\ga_i,$ respectively.
Define $\al: [0,2] \to \CC^2$ 
associated to the 1-form $a$
using (\ref{funk}).
Notice that 
$\int_0^1 \al(t) dt =0$
because the restriction of any element of $\HH^1_A(X; \rrp)$ 
to a loop $\ell:S^1 \to X$ 
is exact whenever $\hol_\ell(A)$ is nontrivial (since 
$\HH^0_{\ell^*(A)}(S^1; \rrp) =0,$
which implies that 
$\HH^1_{\ell^*(A)}(S^1; \rrp) =0$
by Poincar\'e duality). Hence by (\ref{funky})  we see that
$$\int_0^2 \al(t) dt  = \int_1^2 \al(t) dt = L^{-1} \xi_i L.$$
Appealing once again to Corollary \ref{derivative of tr hol} (ii),
it follows that 
\begin{eqnarray*} 
\hess \tr \hol_{\ell\cdot \ga_i}(A)(a,a) 
&=& \int_0^2 \int_0^s \tr \left[ 
L (\al(s) \al(t) + \al(s) \al(t) ) \right] dt ds\\
&=& \int_0^2 \int_0^2 \tr \left[ L \, \al(s) \al(t)\right]  ds dt \\
&=& \tr \left( \int_0^2 L \al(s) ds \int_0^2 \al(t) dt \right) 
= \tr (\xi_i^2 L ) .
\end{eqnarray*}
\end{pf}

Since $a_i=0 \Rightarrow \int_{\ga_i} a = \xi_i =0$, it follows 
from (i) and (ii) above that
the Hessians at $A$ of the real and imaginary parts of
$\tr \hol_{\ga_i}$ and $\tr \hol_{ \ell \cdot \ga_i}$ 
lie in $\UU_{ii}.$
Consider the gauge invariant functions
$f = {\frak Re} \tr \hol_{\ga_i}$ and
$g_\ell = {\frak Im} \tr \hol_{ \ell \cdot \ga_i},$
where ${\frak Re}$ and ${\frak Im}$ denote the real and imaginary parts.
Note that $f$ and $g_\ell$ 
obviously
satisfy condition (i) of Proposition \ref{hessian loops}.
Moreover, since $\hol_{\ga_i}(A)$ is trivial, $D f(A) = 0$.
This follows from formula (i) of Corollary \ref{derivative of tr hol}.
The same formula also implies that the imaginary part
of $D \tr \hol_{ \ell \cdot \ga_i}(A)$ vanishes
since $\tr(L \xi) $ is real for $\xi \in \rrp$ whenever
the $SU(3)$ matrix $L$
is in the image of the standard inclusion $SU(2) \to SU(3)$.
This shows that $D g_\ell(A)=0$,
hence $f$ and $g_\ell$ satisfy 
condition (iii) of Proposition \ref{hessian loops}.
So, we only need to prove that we can span $\UU_{ii}$
with the Hessians of such functions.

For this, we shall use the isomorphism
$V_i \to \rrp$ given by $a \mapsto \int_{\ga_i}a$,
along with the standard identification $\varphi:\rrp \to \CC^2$,
to translate it into a question about symmetric, bilinear pairings
$\CC^2 \times \CC^2 \to \RR.$
Denote by $\langle \cdot , \cdot  \rangle$
the standard complex inner product on $\CC^2$.
If $\xi, \zeta \in \rrp$ and $v, w \in \CC^2$
are given by $v=\varphi(\xi)$ and $w=\varphi(\zeta)$, 
then
$$\tr(\xi \zeta) = - 2 {\frak Re} \langle v,w \rangle.$$
Moreover, if 
$\hat{L} = 
\left( {{\; \al \; \; \; \be}\atop{-\bar{\be} \; \; \bar{\al}}} \right) 
\in SU(2)$
and $L = \hat{L} \oplus 1 \in SU(3),$ then 
$$\tr(L(\xi \zeta + \zeta \xi))=
-\langle \hat{L}(v),w \rangle - \langle \hat{L}(w), v \rangle - 
2 {\frak Re} \langle v,w \rangle.$$
In terms of the real basis 
$\left\{ \left({{1}\atop{0}} \right),
\left({{i}\atop{0}} \right),
\left({{0}\atop{1}} \right),
\left({{0}\atop{i}} \right) \right\}$
for $\CC^2,$
the symmetric bilinear form $\CC^2 \times \CC^2 \to \CC$
given by $$(v,w) \mapsto \langle \hat{L}(v),w \rangle 
+ \langle \hat{L}(w), v \rangle + 2 {\frak Re} \langle v,w \rangle$$
has imaginary part represented by the matrix
\begin{equation} \label{matrix}
\Psi(\hat{L}) = 2 \left[ \begin{array}{rrrr} 
s & 0 & -u & t \\
0 & s & -t & -u \\
-u & -t & -s & 0 \\
t & -u & 0 & -s \\
\end{array} \right].
\end{equation}
where $\al = r+is$ and $\be = t+i u$. 
 
Now $A$ is reducible
(but {\it not} abelian) and thus we have $x,y \in \pi_1(X)$ such that
$\varrho(x)$ and $\varrho(y)$ do not commute.
We claim that the Hessians at $A$ of the four functions
$$f, g_x, g_y, g_{xy}$$
derived from $\ga_{i}$ are linearly independent and form a basis for
the 4-dimensional subspace $\UU_{ii} \subset \UU.$ 

To see this, restrict each Hessian to $V_i \times V_i$
and consider the associated symmetric $4 \times 4$ matrix of the
form (\ref{Uform}).
For example, the matrix associated to
$\hess f(A)$ equals $-2$ times the identity matrix.
Clearly the image of $SU(2)$ under $\Psi$ in (\ref{matrix})
is the complementary subspace of dimension 3.
Thus, it suffices to prove that the Hessians at $A$ of
$g_x, g_y$ and $g_{xy}$,
are linearly independent.  One can see this 
by direct computation;
arranging that  
$\varrho(x)$  is diagonal (by conjugation)
and $\varrho(y)$  is not (by hypothesis), it becomes
a routine exercise in linear algebra. 

This proves that the Hessians at $A$ of $f,g_x,g_y$ and $g_{xy}$
form a basis for $\UU_{ii},$ and to conclude the proof of Proposition \ref{hessian loops},
we need to find, for each $i<j$, functions satisfying (i) and (iii)
whose Hessians span $\UU_{ij}$.

\begin{lem} \label{loopy}
Suppose $i < j$ and $\ell$ is a loop with $\hol_\ell(A)$ nontrivial.
Set $L = \hol_\ell(A) \in SU(3)$.
Suppose further that $a,b \in \HH^1_A(X;\rrp)$,
and set $\xi_k = \int_{\ga_{k}} a \in \rrp$
and 
$\zeta_k = \int_{\ga_k} b \in \rrp$
for $k=1,\ldots, m$.
Then
\begin{itemize}
\item [(i)]  $\hess \tr \hol_{\ga_i \cdot \ga_j} (A)(a, b)= 
\tr((\xi_i+\xi_j)(\zeta_i + \zeta_j))$
\item [(ii)] 
$\hess \tr \hol_{\ga_i \cdot \ell^{-1} \cdot \ga_j\cdot \ell}(A)(a,b)
= \tr( (\xi_i + L \xi_j L^{-1}) ( \zeta_i + L \zeta_j L^{-1} )).$
\item [(iii)] $\hess \tr \hol_{\ell \cdot \ga_i \cdot \ga_j} (A)(a, b)
=\tr(L (\xi_i+\xi_j) ( \zeta_i + \zeta_j))$
\end{itemize}
\end{lem}

\begin{pf} 
By symmetry, it is enough to prove (i)--(iii) 
in the case $a=b.$
For (i), this is just the statement that
$$ \hess \tr \hol_{\ga_i \cdot \ga_j} (A)(a,a) = \tr((\xi_i+\xi_j)^2)$$
for all $a \in \HH^1_A(X;\rrp),$
which follows directly from Corollary
\ref{derivative of tr hol} (ii) as in the proof of Lemma \ref{loopy-0},
using the additional fact that $\int_{\ga_i \cdot \ga_j} a = \xi_i + \xi_j $.

To prove (ii), set
$\be= \ga_i \cdot \ell^{-1} \cdot \ga_j \cdot \ell$
and parameterize $\be$
by the interval $[0,4]$ so that the subintervals
$[0,1], [1,2],[2,3]$ and $[3,4]$ parameterize $\ga_i, \ell^{-1},
\ga_j,$ and $\ell,$
respectively.
Define the function $\al: [0,4] \to \rrp$ 
associated to the 1-form $a$
using (\ref{funk}). 

Now $\al \vert_{[1,2]}$ is exact since 
$\hol_{\ell^{-1}}(A)$ is nontrivial. Similarly, $\al \vert_{[3,4]}$ is exact.
Thus 
$$\int_0^4 \al(t)dt = \int_0^1 \al(t)dt + \int_2^3 \al(t) dt 
=\xi_i + L \xi_j L^{-1}$$
by (\ref{funky}).
Using Corollary \ref{derivative of tr hol} again
and noting that $\hol_\be(A)$ is trivial, it
follows that
\begin{eqnarray*}
 \hess \tr \hol_{\be} (A)(a,a) 
&=&\int_0^4 \int_0^s \tr  (\al(s) \al(t) + \al(t) \al(s) )dt ds \\
&=&\int_0^4 \int_0^4 \tr  \al(s) \al(t) dt ds 
= \tr ((\xi_i + L\xi_j L^{-1})^2).
\end{eqnarray*}
 
To prove part (iii),
set $\be = \ell \cdot \ga_i \cdot \ga_j$
and parameterize
$\be$
by the interval $[0,3]$ so that the subintervals
$[0,1], [1,2]$ and $[2,3]$ parameterize $\ell,\ga_i$ and $\ga_j$, respectively.
Define
the function $\al :[0,3] \to \rrp$
associated to $a$ using (\ref{funk}).
Use (\ref{funky}) and the fact that $\al\vert_{[0,1]}$ is exact 
to conclude that 
$$\int_0^3 \al(t) dt =  
\int_1^3 \al(t) dt = L^{-1} (\xi_i+\xi_j) L.$$
Now Corollary \ref{derivative of tr hol} implies
that
\begin{eqnarray*}
\hess \tr \hol_{\be}(A)(a,a)
&=& \int_0^3 \int_0^3 \tr( L \, \al(s) \al(t)) dt ds \\
&=& \tr \left( \int_0^3   L \, \al(s)  ds \int_0^3 \al(t) dt\right)  
= \tr ((\xi_i+\xi_j)^2 L),
\end{eqnarray*}
and this completes the proof of (iii).
\end{pf}
 
If $a_i =0 = a_j,$ then $\xi_i =0 = \xi_j$ 
and it follows from
(i)--(iii) above that the Hessians at $A$ of the real and imaginary
parts of $\tr \hol_{\ga_i \cdot \ga_j},$
$\tr \hol_{\ell^{-1} \cdot \ga_i \cdot \ell \cdot \ga_j}$,
and $\tr \hol_{\ell \cdot \ga_i \cdot \ga_j}$ lie in $\UU_{ij}.$
Consider the gauge invariant functions
$\AA \to \RR$ defined by
$f = {\frak Re} \tr \hol_{\ga_i \cdot \ga_j},$ 
$f_\ell= {\frak Re} \tr \hol_{\ell^{-1} \cdot \ga_i \cdot \ell \cdot \ga_j}$ 
and $g_\ell= {\frak Im} \tr \hol_{\ell \cdot \ga_i \cdot \ga_j}.$
Then conditions (i) and (iii) of Proposition \ref{hessian loops}
are satisfied for $f, f_\ell$ and $g_\ell.$
Condition (i) obviously holds,
and condition (iii)
follows from Corollary \ref{derivative of tr hol} 
just as in the case $i=j$ since the loops for $f$ and $f_\ell$
(coming from parts (i) and (ii) 
of Lemma \ref{loopy}) have trivial holonomy
and since $g_\ell$ is the {\it imaginary} part of trace of holonomy.

So, to complete the proof of
\ref{hessian loops},
we just need to
show that we can span $\UU_{ij}$ with the
Hessians of such functions.
Restricting 
elements in $\UU_{ij}$ to $V_i \times V_j$
we obtain $4 \times 4$ matrices as in (\ref{Uform}).
In contrast to the previous case when $i=j,$
these matrices are not generally symmetric.

Suppose $a \in V_i$ and $b \in V_j$. 
Then $\xi_j =0$ and $\zeta_i =0$.
Let $v =\varphi(\xi_{i}) \in \CC^{2}$ and
$w = \varphi(\zeta_{j}) \in \CC^{2}.$
If $\hat{L} = \hol_{\ell}¥(\hat{A}) \in SU(2)$,
(so $L = \hat{L} \oplus 1$), then Lemma \ref{loopy} implies that
\begin{eqnarray*}
\hess f(A)(a,b) &=& \tr(\xi_i \zeta_j) = -2 {\frak Re} \langle v,w \rangle, \\
\hess f_\ell(A)(a,b) &=& \tr( \xi_i \, L \, \zeta_j L^{-1})=
- 2 {\frak Re} \langle v, \hat{L}(w)\rangle, \\
\hess g_\ell(A)(a,b) &=& {\frak Im} \tr(L\, \xi_i \zeta_j)=
-{\frak Im}(\langle \hat{L}(v), w \rangle + \langle w,v \rangle).
\end{eqnarray*}
Writing 
$\hat{L} = 
\left( {{ \; \al \; \; \;  \be} \atop {-\bar{\be} \; \; \bar{\al}}} \right)$
where $\al = r+is$ and $\be = t+iu,$
then in terms of the real basis
$\left\{ \left({{1}\atop{0}} \right),
\left({{i}\atop{0}} \right),
\left({{0}\atop{1}} \right),
\left({{0}\atop{i}} \right) \right\}$
for $\CC^2,$
the bilinear pairing $\CC^2 \times \CC^2 \to \CC$
given by 
$(v,w) \mapsto \langle v, \hat{L}(w) \rangle$
has real part represented by the matrix
\begin{equation} \label{matrix2}
\Phi(\hat{L}) = \left[\begin{array}{cccc}
r & -s &  t & -u \\
s & r &  u &  t \\
-t & -u & r & s \\
u & -t & -s & r
\end{array}\right].
\end{equation}
Likewise, the bilinear pairing $\CC^2 \times \CC^2 \to \CC$
given by 
$(v,w) \mapsto \langle \hat{L}(v), w \rangle + \langle w,v \rangle$
has imaginary part represented by the matrix
\begin{equation} \label{matrix3}
\Psi(\hat{L}) =\left[\begin{array}{cccc}
s & 1-r & u &  t \\
r-1 & s & -t & u \\
u & -t & -s & 1-r \\
t & u & r-1 & -s 
\end{array} \right].
\end{equation}
Notice that the images of $SU(2)$ under $\Phi$ and $\Psi$
span complementary 4-dimensional subspaces of the 8-dimensional
space of matrices
of the form (\ref{Uform}).

Choose $x,y \in \pi_1(X)$ as before
so that $\varrho(x)$ and $\varrho(y)$
do not commute.
We first claim that the Hessians at $A$ of $f, f_x, f_y$ and $f_{xy}$
are linearly independent. In fact, after restricting to
$V_i \times V_j,$ they span the 4-dimensional subspace
of matrices (\ref{matrix2}).
To show this, one only needs to show that the image of the set
$\{ I, \varrho(x), \varrho(y), \varrho(xy) \}$
under $\Phi$ is linearly independent. Again, this 
follows from the hypotheses on $\varrho(x)$
and $\varrho(y)$ easily
after assuming (by conjugation) that
$\varrho(x)$ is diagonal.

The complementary 4-dimensional subspace of $\UU_{ij}$
given by (\ref{matrix3}) can be spanned using functions $g_\ell$.
The image of the set
$\{ I, \varrho(x), \varrho(y), \varrho(xy) \}$ 
under $\Psi$ is  linearly  dependent because
$\Psi(I) =0.$ However, a straightforward check shows that
the image of $\{ \varrho(x), \varrho(x^2), \varrho(y), \varrho(xy) \}$
under $\Psi$ is linearly independent.
Hence, it follows that the Hessians of $g_x, g_{x^2}, g_{y}$ and $g_{x y}$
are linearly independent. Since their span is complementary to that
of the Hessians at $A$ of $f, f_x, f_y$ and $f_{xy}$,
together they span $\UU_{ij}$
and this completes the proof of Proposition \ref{hessian loops}.
\end{pf}

\end{document}